\documentclass[11pt,a4paper,reqno]{amsart}

\usepackage[utf8]{inputenc}
\usepackage[T1]{fontenc}
\usepackage[activeacute,english]{babel}
\usepackage{lmodern}
\usepackage{microtype}

\usepackage{amsmath,amssymb,amsthm,amsfonts,mathtools}

\usepackage{enumitem}
\usepackage[colorlinks,bookmarks,linkcolor=black,citecolor=black]{hyperref}
\usepackage{fullpage}

\setenumerate{leftmargin=*}

\newtheorem{theorem}{Theorem}[section]
\newtheorem{lemma}[theorem]{Lemma}
\newtheorem{corollary}[theorem]{Corollary}
\newtheorem{claim}[theorem]{Claim}
\newtheorem{conjecture}{Conjecture}
\theoremstyle{definition}
\newtheorem{definition}[theorem]{Definition}
\newtheorem{setup}[theorem]{Setup}

\DeclareMathOperator{\conn}{conn}
\DeclareMathOperator{\red}{red}
\DeclareMathOperator{\blue}{blue}

\newcommand{\bQ}{\mathbf{Q}}
\newcommand{\bd}{\mathbf{d}}
\newcommand{\cC}{\mathcal{C}}
\newcommand{\cH}{\mathcal{H}}
\newcommand{\cJ}{\mathcal{J}}
\newcommand{\cM}{\mathcal{M}}
\newcommand{\cP}{\mathcal{P}}
\newcommand{\Jchat}{\hat{\cJ}}
\newcommand{\NN}{\mathbb{N}}
\newcommand{\eps}{\varepsilon}
\newcommand{\nth}{\varnothing}
\newcommand{\sm}{\setminus}

\def\endofClaim{\hfill\scalebox{.6}{$\blacksquare$}}
\newcommand{\oldqed}{}
\newenvironment{claimproof}[1][Proof]{
  \renewcommand{\oldqed}{\qedsymbol}
  \renewcommand{\qedsymbol}{\endofClaim}
  \begin{proof}[#1]
}{
  \end{proof}
  \renewcommand{\qedsymbol}{\oldqed}
}

\title{Off-diagonal Ramsey numbers for $3$-uniform tight hypergraph cycles}
\author{Eng Keat Hng}
\address{Extremal Combinatorics and Probability Group (ECOPRO), Institute for Basic Science (IBS), Daejeon, South Korea.}
\email{hng@ibs.re.kr}
\thanks{Research supported by the Institute for Basic Science (IBS-R029-C4).}
\date{}

\begin{document}

\begin{abstract}
    Let~$C^{(3)}_n$ be the~\emph{$3$-uniform tight cycle on~$n$ vertices}, that is, the~$3$-uniform hypergraph on a cyclically ordered set of~$n$ vertices whose edges are precisely the triples of consecutive vertices in the given cyclic order. We asymptotically determine the off-diagonal Ramsey numbers of~$3$-uniform tight cycles. More specifically, we show that for~$s\ge t\ge13$ and~$i,j\in\{0,1,2\}$ we have
    \begin{equation*}
        r(C^{(3)}_{3s+i},C^{(3)}_{3t+j})=
        \begin{cases}
            3s+t+o(s) &\textrm{ if }i=0,\ j=0 \\
            \max\{5s,3(s+t)\}+o(s) &\textrm{ if }i\ne0,\ j\ne0 \\
            5s+o(s) &\textrm{ if }i=0,\ j\ne0 \\
            \max\{3s+t,5t\}+o(s) &\textrm{ if }i\ne0,\ j=0\,.
        \end{cases}
    \end{equation*}
    This is an off-diagonal analogue of a result of Haxell, {\L}uczak, Peng, R\"{o}dl, Ruci\'{n}ski and Skokan [Combin.\ Probab.\ Comput.\ 18 (2009), no.\ 1-2, 165--203] on the diagonal Ramsey numbers of~$3$-uniform tight cycles.
\end{abstract}

\maketitle

\section{Introduction}

For~$k$-uniform hypergraphs~$H$ and~$H'$ the \emph{Ramsey number}~$r(H,H')$ is the least positive integer~$N$ such that every red-blue edge-colouring of the complete~$k$-uniform hypergraph~$K^{(k)}_N$ contains a red copy of~$H$ or a blue copy of~$H'$. The study of cycle Ramsey numbers is a major research direction in Ramsey theory, and a recurring theme in this direction is that the leading term behaviour of cycle Ramsey numbers and the associated extremal colourings often exhibit a strong dependency on simple arithmetic properties. In graphs, the relevant property is parity; indeed, a classical theorem~\cite{BondyErdos,FaudreeSchelp,Rosta} states that for all~$m\ge n\ge 3$ we have
\begin{equation*}
    r(C_m,C_n)=
    \begin{cases}
        2m-1 & \textrm{ if $n$ is odd and }(m,n)\ne(3,3) \\
        m-1+n/2 & \textrm{ if $m$ and $n$ are even and }(m,n)\ne(4,4) \\
        \max(m-1+n/2,2n-1) & \textrm{ if $m$ is odd and $n$ is even.}
    \end{cases}
\end{equation*}
A similar phenomenon is also observed with more colours. Indeed, for the~$q$-colour Ramsey number with~$q\ge3$, Knierim and Su~\cite{KnierimSu} proved that~$r_q(C_n)\le\left(q-\frac{1}{2}+o(1)\right)n$ for all even~$n$, while Jenssen and Skokan~\cite{JenssenSkokan} proved that~$r_q(C_n)=2^{q-1}(n-1)+1$ for all sufficiently large odd~$n$.

Ramsey numbers have also been studied in the context of~$k$-uniform hypergraphs, where a variety of different notions of cycles coexist. Here we will focus on \emph{tight cycles} and refer the reader to~\cite{Shahsiah} and the references therein for Ramsey numbers of loose cycles. Let~$C^{(k)}_n$ denote the~\emph{$k$-uniform tight cycle on~$n$ vertices}, that is, the~$k$-uniform hypergraph on a cyclically ordered set of~$n$ vertices whose edges are precisely the sets of~$k$ consecutive vertices in the given cyclic order. Haxell, {\L}uczak, Peng, R\"{o}dl, Ruci\'{n}ski and Skokan~\cite{HaxellLuczakPengRodl} asymptotically determined the diagonal Ramsey numbers of~$3$-uniform tight cycles. More precisely, they showed that
\begin{equation*}
    r(C^{(3)}_{3n+i},C^{(3)}_{3n+i}) = 
    \begin{cases}
        (4+o(1))n & \textrm{ if }i=0 \\
        (6+o(1))n & \textrm{ if }i\in\{1,2\}\,.
    \end{cases}
\end{equation*}
They also proposed the following general conjecture for all~$k\ge3$ about the asymptotics and arithmetic features of~$r(C^{(k)}_n,C^{(k)}_n)$. Since every pair of consecutive edges in a~$k$-uniform tight cycle intersects in exactly~$k-1$ vertices, it is natural that the residue class~$i$ modulo~$k$ of the tight cycle length features prominently in the relevant arithmetic property~$\gcd(k,i)$.

\begin{conjecture}[Haxell, {\L}uczak, Peng, R\"{o}dl, Ruci\'{n}ski and Skokan~\cite{HaxellLuczakPengRodl}] \label{conj:k-cycle-ramsey}
    For all~$k\ge3$ and~$0\le i<k$, writing~$d:=\gcd(k,i)$, we have
    \begin{equation*}
        r(C^{(k)}_{kn+i},C^{(k)}_{kn+i}) = (1+o(1))\frac{d+1}{d}kn\,.
    \end{equation*}
\end{conjecture}

Recently, Lo and Pfenninger~\cite{LoPfenninger2025} proved the divisible case~$i=0$ of Conjecture~\ref{conj:k-cycle-ramsey} for~$k=4$, and Pfenninger~\cite{Pfenninger2025} subsequently extended it to general~$k\ge3$. The non-divisible cases~$i\ne0$ of Conjecture~\ref{conj:k-cycle-ramsey} remain widely open; the central source of difficulty is arithmetic in nature.

In this paper, we turn our attention to off-diagonal Ramsey numbers for tight cycles. We prove the following theorem, which asymptotically determines the off-diagonal Ramsey numbers of $3$-uniform tight cycles and builds on the work in~\cite{HaxellLuczakPengRodl}. The off-diagonal picture features highly non-trivial interaction between the two residue classes and the off-diagonal ratio~$t/s$ with highly asymmetric mixed residue regimes and different ranges dominated by different lower-bound colourings.

\begin{theorem} \label{thm:main}
    The following hold.
    \begin{enumerate} [label=(\alph{*})]
        \item \label{item:thm-lower-bound} For all $s\ge t\ge 1$ and $i,j\in\{0,1,2\}$ we have
        \begin{equation*}
            r(C^{(3)}_{3s+i},C^{(3)}_{3t+j})\ge
            \begin{cases}
                3s+t-1 &\textrm{ if }i=0,\ j=0 \\
                \max(5s+2i-1,3(s+t)+i+j-1) &\textrm{ if }i\ne0,\ j\ne0 \\
                5s-1 &\textrm{ if }i=0,\ j\ne0 \\
                \max(3s+t+i-1,5t-1) &\textrm{ if }i\ne0,\ j=0\,.
            \end{cases}
        \end{equation*}
        \item \label{item:thm-upper-bound} For all~$\eta>0$ there is~$s_0\in\NN$ such that for all~$s\ge s_0$, $13\le t\le s$ and~$i,j\in\{0,1,2\}$ we have
        \begin{equation*}
            r(C^{(3)}_{3s+i},C^{(3)}_{3t+j})\le
            \begin{cases}
                3s+t + \eta s &\textrm{ if }i=0,\ j=0 \\
                \max(5s,3(s+t)) + \eta s &\textrm{ if }i\ne0,\ j\ne0 \\
                5s + \eta s &\textrm{ if }i=0,\ j\ne0 \\
                \max(3s+t,5t) + \eta s &\textrm{ if }i\ne0,\ j=0\,.
            \end{cases}
        \end{equation*}
    \end{enumerate}
\end{theorem}

The four cases in Theorem~\ref{thm:main} exhibit distinct structural behaviour. When both lengths are divisible by~$3$, the leading term is~$3s+t$; this strongly resembles the even--even case in graphs. When neither length is divisible by~$3$, the two obstructions~$5s$ and~$3(s+t)$ compete and produce a smooth transition as~$t/s$ varies. If only the blue length is non-divisible, then the~$5s$ obstruction dominates throughout; in the other mixed case, the~$3s+t$ obstruction competes with the~$5t$ obstruction. In particular, Theorem~\ref{thm:main} completely describes the off-diagonal asymptotic behaviour and recovers both diagonal constants~$4$ and~$6$.

The lower bounds in Theorem~\ref{thm:main}\ref{item:thm-lower-bound} are witnessed by three explicitly defined colourings; the first two are straightforward extensions of the constructions given in~\cite{HaxellLuczakPengRodl}, while the third colouring requires an additional vertex class which features most prominently in the mixed residue regimes. The upper bound applies the well-established connected matching method of {\L}uczak, which reduces the problem of finding monochromatic tight cycles of divisible lengths to that of finding a sufficiently large matching in a monochromatic component; for non-divisible lengths we additionally require a residue-adjusting gadget (in our setting, a copy of~$C^{(3)}_4$ or~$C^{(3)}_5$) in the same component. The main difficulty in our proof lies in forcing these gadgets at the desired scale of~$5s$ instead of the~$6s$-type bound one would obtain by directly adapting the arguments in~\cite{HaxellLuczakPengRodl}. In particular, the proof builds on the methods in~\cite{HaxellLuczakPengRodl} by developing innovative matching-pivot techniques and monochromatic component structural analysis.

As a straightforward corollary of Theorem~\ref{thm:main}, we also asymptotically determine the off-diagonal Ramsey numbers for~$3$-uniform tight paths; let~$P^{(k)}_n$ denote the~\emph{$k$-uniform tight path on~$n$ vertices}, that is, the~$k$-uniform hypergraph on a linearly ordered set of~$n$ vertices whose edges are precisely the sets of~$k$ consecutive vertices in the given linear order.

\begin{corollary} \label{cor:paths}
    The following hold.
    \begin{enumerate} [label=(\alph{*})]
        \item \label{item:cor-lower-bound} For all $s\ge t\ge 3$ we have~$r(P^{(3)}_s,P^{(3)}_t) \ge s+\left\lfloor\frac{t}{3}\right\rfloor-1$.
        \item \label{item:cor-upper-bound} For all~$\eta>0$ there is~$s_0\in\NN$ such that for all~$s\ge s_0$ and~$39\le t\le s$ we have~$r(P^{(3)}_s,P^{(3)}_t) \le (1+\eta)s+\left\lfloor\frac{t}{3}\right\rfloor$.
    \end{enumerate}
\end{corollary}

Indeed, Corollary~\ref{cor:paths}\ref{item:cor-upper-bound} follows immediately from the~$i=j=0$ case of Theorem~\ref{thm:main}\ref{item:thm-upper-bound} and Corollary~\ref{cor:paths}\ref{item:cor-lower-bound} follows by the same proof as that for Theorem~\ref{thm:main}\ref{item:thm-lower-bound} in Section~\ref{ssec:lower-bounds}.

The rest of the paper is organised as follows. In Section~\ref{ssec:lower-bounds} we give the extremal colourings and prove the lower bounds, while in Section~\ref{ssec:upper-bounds} we outline the hypergraph version of the connected matching method of {\L}uczak which we use. In Section~\ref{sec:matching-xdivs} we introduce the structural analysis from~\cite{HaxellLuczakPengRodl}, develop the key techniques needed to force the necessary gadgets in the non-divisible cases and apply them together with the basic connected matching lemma (Lemma~\ref{lem:matching-div-div}) to prove the connected matching lemmas with gadgets (Lemmas~\ref{lem:matching-xdiv-xdiv}--\ref{lem:matching-xdiv-div}). In Section~\ref{sec:regularity-cycle-embedding} we provide a proof of the technical lemma (Lemma~\ref{lem:monox-cycle-embedding}) which encapsulates the problem reduction via the hypergraph connected matching method. In Appendix~\ref{sec:matching-div-div} we provide a proof of Lemma~\ref{lem:matching-div-div}, and in Appendix~\ref{sec:matching-xdivs} we provide the deferred proofs of a variety of lemmas in Section~\ref{sec:matching-xdivs}.

\subsection{Notation and definitions}

Write~$\NN$ for the set of positive integers and~$\NN_0$ for the set~$\NN\cup\{0\}$. For~$m\in\NN_0$ write~$[m]$ for the set~$\{1,\dots,m\}$ and~$[m]_0$ for the set~$[m]\cup\{0\}$. For a set~$S$ and~$m\in\NN_0$ let~$\binom{S}{m}$ denote the set of subsets of~$S$ of size~$m$. For positive real numbers~$x$ and~$y$, write~$x\pm y$ to mean an appropriate real number~$z$ satisfying~$x-y\le z\le x+y$. For example, for positive real numbers~$w$, $x$, $y$ and~$a$ the equation~$w=(x\pm y)a$ would mean~$(x-y)a\le w\le(x+y)a$. We will use constant hierarchy notation, that is, when we say that a statement holds for~$0<a \ll b \le 1$ it means that there is a non-decreasing function~$a_0:(0,1]\to(0,1]$ such that the statement holds for all~$a,b\in(0,1]$ such that~$a\le a_0(b)$; hierarchies with more constants are defined analogously, and we implicitly assume that if~$\frac{1}{m}$ is a term in a constant hierarchy then~$m$ is a positive integer. We omit floors and ceilings whenever they do not affect the calculations.

A hypergraph~$\cH$ consists of a vertex set~$V(\cH)$ and an edge set~$E(\cH)$ whose elements are subsets of~$V(\cH)$. We often identify a hypergraph~$\cH$ with its edge set~$E(\cH)$. A~\emph{$k$-uniform hypergraph}, or~\emph{$k$-graph}, is a hypergraph in which every edge has size~$k$. Let~$K^{(k)}_n$ denote the complete~$k$-graph on~$n$ vertices. A \emph{matching} of \emph{size}~$s$ in a~$k$-graph~$H$ is a set of~$s$ pairwise disjoint edges of~$H$. For a matching~$M$ we write~$V(M)$ for the set of vertices covered by~$M$.

\section{Lower bounds and proof outline for upper bounds}

The main work in proving Theorem~\ref{thm:main} lies in the upper bounds. In this section, we first prove the lower bounds in Section~\ref{ssec:lower-bounds} and then introduce the main ingredients and lemmas needed to prove the upper bounds in Section~\ref{ssec:upper-bounds}. We will apply the main lemmas to give a proof of the upper bounds at the end of this section.

\subsection{Proving lower bounds} \label{ssec:lower-bounds}

In this subsection, we prove Theorem~\ref{thm:main}\ref{item:thm-lower-bound}, which gives lower bounds on~$r(C^{(3)}_{3s+i},C^{(3)}_{3t+j})$. The proof involves giving three constructions of~$2$-edge-coloured~$3$-graphs; the first two are straightforward generalisations of the constructions given in~\cite{HaxellLuczakPengRodl} and the third one involves a notable adjustment to the second construction.

\begin{proof}[Proof of Theorem~\ref{thm:main}\ref{item:thm-lower-bound}]
We shall prove several sets of lower bounds which collectively give the desired lower bounds; each set arises from a construction and applies to specific values of~$i$ and~$j$.

To show that~$r(C^{(3)}_{3s+i},C^{(3)}_{3t}) \ge 3s+t+i-1$ for~$i\in\{0,1,2\}$, let~$K$ be a copy of~$K^{(3)}_{3s+t+i-2}$, fix a partition~$V(K)=X\sqcup Y$ with~$|X|=3s+i-1$ and~$|Y|=t-1$, and colour the edges of~$K$ such that an edge~$e\in E(K)$ is red if~$e\subseteq X$ and blue otherwise. Suppose that this edge colouring contains a red copy~$C$ of~$C^{(3)}_{3s+i}$. Since~$|X|<3s+i$, there is an edge of~$C$ which contains a vertex of~$Y$. But this edge must be coloured blue, yielding a contradiction. Now suppose that this edge colouring contains a blue copy of~$C^{(3)}_{3t}$. Then it contains~$t$ disjoint blue edges. Each blue edge must contain a distinct vertex of~$Y$, so~$Y$ contains at least~$t$ vertices, which is a contradiction.

To show that~$r(C^{(3)}_{3s+i},C^{(3)}_{3t+j}) \ge 3(s+t)+i+j-1$ for~$i,j\in\{1,2\}$, let~$K$ be a copy of~$K^{(3)}_{3(s+t)+i+j-2}$, fix a partition~$V(K)=X\sqcup Y$ with~$|X|=3s+i-1$ and~$|Y|=3t+j-1$, and colour an edge~$e\in E(K)$ red if~$|e\cap X|$ is odd and blue otherwise. Suppose that this edge colouring contains a red copy~$C$ of~$C^{(3)}_{3s+i}$ with cyclic vertex sequence~$v_1,\dots,v_{\ell}$, and for all~$h\in[\ell]$ set~$a_h=1$ when~$v_h\in X$ and~$a_h=0$ otherwise. Since~$C$ is monochromatic red, it follows that the parity of the sum of any three (cyclically) consecutive terms must be the same, and so for all~$h\in[3s+i]$ we have~$a_{h+3}=a_h$ (with the indices taken modulo~$\ell$). Now~$\gcd(3,3s+i)=1$ and any red edge~$e$ has odd~$|e\cap X|$, so we have~$a_1=\dots=a_{3s+i}=1$. Hence, $X$ contains all~$3s+i$ vertices of~$C$, which yields a contradiction. Finally, an analogous argument implies that the given edge colouring also contains no blue copy of~$C^{(3)}_{3t+j}$.

To show that $r(C^{(3)}_{3s+i},C^{(3)}_{3t+j}) \ge 5s+2i-1$ for~$i\in\{0,1,2\}$ and~$j\in\{1,2\}$, let~$K$ be a copy of~$K^{(3)}_{5s+2i-2}$ and fix a partition $V(K)=X\sqcup Y\sqcup Z$ with $|X|=3s+i-1$, $|Y|=\min(3t+j,2s+i)-1$ and~$|Z|=\max(0,2s+i-3t-j)$. Hence, we have~$|Y\cup Z|=2s+i-1$. Colour an edge~$e\in E(K)$ blue if~$|e\cap X|=2$ or~$e\subseteq Y$, and red otherwise. Suppose that this edge colouring contains a red copy~$C$ of~$C^{(3)}_{3s+i}$ with cyclic vertex sequence~$v_1,\dots,v_{3s+i}$, and for all~$h\in[3s+i]$ set~$a_h=1$ when~$v_h\in X$ and~$a_h=0$ otherwise. Since~$C$ is monochromatic red and~$|X|<3s+i$, the cyclic sequence~$(a_h)_{h\in[3s+i]}$ contains at least one~$0$ and has no~$3$-term interval with exactly two~$1$s. It follows that~any two~$1$s in~$(a_h)_{h\in[3s+i]}$ is separated by at least two~$0$s, so~$C$ contains at least~$\left\lceil\frac{2(3s+i)}{3}\right\rceil=2s+i$ vertices outside~$X$, which contradicts~$|Y\cup Z|=2s+i-1$. Now suppose that the edge colouring contains a blue copy~$C$ of~$C^{(3)}_{3t+j}$ with cyclic vertex sequence~$v_1,\dots,v_{3t+j}$, and for all~$h\in[3t+j]$ set~$a_h=1$ when~$v_h\in X$ and~$a_h=0$ otherwise. Suppose further that~$C$ has an edge entirely contained in~$Y$. Then it follows by a vertex-by-vertex argument around~$C$ that~$Y$ contains all~$3t+j$ vertices of~$C$, which contradicts~$|Y|<3t+j$. Hence, every edge of~$C$ has exactly two vertices in~$X$, so for all~$h\in[3t+j]$ we have~$a_{h+3}=a_h$ (with the indices taken modulo~$\ell$). Since~$\gcd(3,3t+j)=1$, we have~$a_1=\dots=a_{3t+j}$. Hence, every edge of~$C$ has exactly zero or three vertices in~$X$, which is a contradiction. Finally, by interchanging the two colours and the pairs~$(s,i)$ and~$(t,j)$, we obtain~$r(C^{(3)}_{3s+i},C^{(3)}_{3t+j}) \ge 5t+2j-1$ for~$i\in\{1,2\}$ and~$j\in\{0,1,2\}$.
\end{proof}

\subsection{Proving upper bounds} \label{ssec:upper-bounds}

Our proof of Theorem~\ref{thm:main}\ref{item:thm-upper-bound} uses a hypergraph version of the well-established technique often known as the \emph{connected matching method of {\L}uczak}~\cite{Luczak}. Broadly speaking, the idea is to apply a suitable regularity lemma to the given (hyper)graph, obtain a suitably connected matching (possibly with some additional gadgets) in the so-called reduced graph and convert this matching in the reduced graph into the desired cycle in the original (hyper)graph. This technique reduces the original Ramsey problem of obtaining monochromatic tight cycles in a~$2$-edge-coloured~$3$-graph to a Ramsey-type problem about obtaining suitably sized monochromatic matchings (possibly with some additional gadgets) which are appropriately connected.

We now introduce the appropriate notion of connectivity. Let~$H$ be a~$k$-graph. We say that two edges~$e$ and~$e'$ of~$H$ are \emph{tightly connected} if there exists a sequence~$e_0,\dots,e_{\ell}$ of edges of~$H$ such that~$e_0=e$, $e_{\ell}=e'$ and for all~$i\in[\ell]$ we have~$|e_{i-1}\cap e_{i}|=k-1$. It is straightforward to check that the property of being connected defines an equivalence relation on the set of edges of~$H$, and we call the equivalence classes the \emph{tight components} of $H$; since the notion of connectedness defined above is the only one we work with, we shall often omit the modifiers `tight' and `tightly' throughout this paper. Observe that for~$k=2$ this corresponds to the usual notion of connectedness in graphs. A \emph{connected matching} in~$H$ is a matching whose edges are all contained in a single component of~$H$.

A~\emph{$2$-edge-coloured~$k$-graph} is a~$k$-graph together with a red-blue colouring of its edges. Let~$H$ be a~$2$-edge-coloured~$k$-graph. Write~$H_{\red}$ and~$H_{\blue}$ for the~$k$-graphs on~$V(H)$ induced by the red and blue edges of~$H$, respectively. A subgraph of~$H$ is \emph{monochromatic} if it is a subgraph of~$H_{\red}$ or~$H_{\blue}$. A \emph{red component (resp.\ blue component)} of~$H$ is a component of~$H_{\red}$ (resp.~$H_{\blue}$); a \emph{monochromatic component} of~$H$ is a red component or blue component of~$H$.

Let $H$ be a $3$-graph. Write $\partial H$ for the set of all pairs $xy$ for which there is $z$ such that $xyz\in E(H)$; we say that each~$xy\in\partial H$ is \emph{active}. Observe that for any two distinct components $H_1$ and $H_2$ of $H$ we have
\begin{equation} \label{eq:components-disjoint-shadows}
    \partial H_1\cap\partial H_2=\nth\,.
\end{equation}
For a vertex~$x\in V(H)$ set~$N_{\partial H}(x):=\{y\in V(H):xy\in\partial H\}$, and for vertices~$x,y\in V(H)$ set~$N_H(x,y):=\{z\in V(H):xyz\in E(H)\}$. We sometimes omit the subscripts~$\partial H$ and~$H$ when they are clear from the context. We say that a~$3$-graph~$H$ on~$N$ vertices is~\emph{$\eps$-complete} if for each~$x\in V(H)$ we have~$N_{\partial H}(x)\ne\nth$ and for each~$xy\in\partial H$ we have
\begin{equation} \label{eq:complete-codegree}
    |N_H(x,y)|\ge(1-\eps)N\,.
\end{equation}

The key ingredient in the proof of Theorem~\ref{thm:main}\ref{item:thm-upper-bound} is the following series of four lemmas, which resolve the reduced problem and correspond to the four cases in the statement of Theorem~\ref{thm:main}\ref{item:thm-upper-bound}. Let~$M^{(3)}_{s,\conn}$ denote a connected matching of size~$s$ and~$M^{(3)}_{s,\conn,+}$ denote a connected matching of size~$s$ with a copy of~$C^{(3)}_4$ or~$C^{(3)}_5$ in the same monochromatic component. The basic connected matching lemma, Lemma~\ref{lem:matching-div-div}, is an easy generalisation of~\cite[Lemma~2.1]{HaxellLuczakPengRodl}, and for completeness we provide a proof of this lemma in Appendix~\ref{sec:matching-div-div}. Indeed, the main work is in developing the techniques needed to obtain the~$5s$-type bounds seen in Lemmas~\ref{lem:matching-xdiv-xdiv}--\ref{lem:matching-xdiv-div}. We provide the proofs of these lemmas in Section~\ref{sec:matching-xdivs}.

\begin{lemma} \label{lem:matching-div-div}
    Let~$s,t\in\NN$ and~$\eps,\eta>0$ satisfy~$\frac{1}{s}\ll\eps\ll\eta$ and~$s\ge t$. Suppose that~$H$ is a~$2$-edge-coloured~$\eps$-complete~$3$-graph on~$N=3s+t+\eta s$ vertices. Then $H$ contains a red~$M^{(3)}_{s,\conn}$ or a blue~$M^{(3)}_{t,\conn}$.
\end{lemma}

\begin{lemma} \label{lem:matching-xdiv-xdiv}
    Let~$s,t\in\NN$ and~$\eps,\eta>0$ satisfy~$\frac{1}{s}\ll\eps\ll\eta$ and~$s\ge t$. Suppose that~$H$ is a~$2$-edge-coloured~$\eps$-complete~$3$-graph on~$N=\max(5s,3(s+t))+\eta s$ vertices. Then $H$ contains a red~$M^{(3)}_{s,\conn,+}$ or a blue~$M^{(3)}_{t,\conn,+}$.
\end{lemma}

\begin{lemma} \label{lem:matching-div-xdiv}
    Let~$s,t\in\NN$ and~$\eps,\eta>0$ satisfy~$\frac{1}{s}\ll\eps\ll\eta$ and~$s\ge t$. Suppose that~$H$ is a~$2$-edge-coloured~$\eps$-complete~$3$-graph on~$N=5s+\eta s$ vertices. Then $H$ contains a red~$M^{(3)}_{s,\conn}$ or a blue~$M^{(3)}_{t,\conn,+}$.
\end{lemma}

\begin{lemma} \label{lem:matching-xdiv-div}
    Let~$s,t\in\NN$ and~$\eps,\eta>0$ satisfy~$\frac{1}{s}\ll\eps\ll\eta$ and~$s\ge t$. Suppose that~$H$ is a~$2$-edge-coloured~$\eps$-complete~$3$-graph on~$N=\max(3s+t,5t)+\eta s$ vertices. Then $H$ contains a red~$M^{(3)}_{s,\conn,+}$ or a blue~$M^{(3)}_{t,\conn}$.
\end{lemma}

We need the following definition to state the hypergraph version of the connected matching method of {\L}uczak which we use.

\begin{definition}
    Let~$m^*_1(\beta,\eps,n)$ be the largest~$m$ such that every~$2$-edge-coloured~$\eps$-complete~$3$-graph on~$n$ vertices contains a matching of size at least~$m$ whose edges all lie in a single red tight component or a matching of size at least~$\beta m$ whose edges all lie in a single blue tight component. Let~$m^*_2(\beta,\eps,n)$ be the largest~$m$ such that every~$2$-edge-coloured~$\eps$-complete~$3$-graph on~$n$ vertices contains a matching of size at least~$m$ whose edges all lie in a single red tight component which contains a copy of~$C^{(3)}_4$ or~$C^{(3)}_5$, or a matching of size at least~$\beta m$ whose edges all lie in a single blue tight component which contains a copy of~$C^{(3)}_4$ or~$C^{(3)}_5$. Let~$m^*_3(\beta,\eps,n)$ be the largest~$m$ such that every~$2$-edge-coloured~$\eps$-complete~$3$-graph on~$n$ vertices contains a matching of size at least~$m$ whose edges all lie in a single red tight component or a matching of size at least~$\beta m$ whose edges all lie in a single blue tight component which contains a copy of~$C^{(3)}_4$ or~$C^{(3)}_5$. Let~$m^*_4(\beta,\eps,n)$ be the largest~$m$ such that every~$2$-edge-coloured~$\eps$-complete~$3$-graph on~$n$ vertices contains a matching of size at least~$m$ whose edges all lie in a single red tight component which contains a copy of~$C^{(3)}_4$ or~$C^{(3)}_5$, or a matching of size at least~$\beta m$ whose edges all lie in a single blue tight component. Given~$i\in[4]$ and~$f\colon S\to(0,1]$ with~$S\subseteq(0,1]$, set
    \[m^*_i(f):=\liminf_{\eps\to0}\liminf_{n\to\infty}\inf_{\beta\in S}\frac{m^*_i(\beta,\eps,n)}{f(\beta)n}\,.\]
    For~$\beta\in(0,1]$ set~$\mu^*_1(\beta):=\frac{1}{3+\beta}$, $\mu^*_2(\beta):=\frac{1}{\max(5,3(1+\beta))}$, $\mu^*_3(\beta):=\frac{1}{5}$ and~$\mu^*_4(\beta):=\frac{1}{\max(5\beta,3+\beta)}$.
\end{definition}

The following lemma encapsulates the problem reduction from obtaining monochromatic tight cycles in a~$2$-edge-coloured complete~$3$-graph to obtaining tightly connected matchings in a~$2$-edge-coloured~$\eps$-complete~$3$-graph. We provide its proof in Section~\ref{sec:regularity-cycle-embedding}.

\begin{lemma} \label{lem:monox-cycle-embedding}
    Let~$\delta>0$ and~$n\in\NN$ satisfy~$\frac{1}{n}\ll\delta\ll1$. Suppose that~$K=K_{\red}\cup K_{\blue}$ is a~$2$-edge-coloured complete~$3$-graph on~$n$ vertices. Then the following hold for all functions~$f\colon S\to(0,1]$ with~$S\subseteq(0,1]$ and~$\beta\in S$.
    \begin{enumerate}[label=(\roman{*})]
        \item \label{item:monox-cycle-embedding-div-div} $K$ contains a red tight cycle of length~$\ell$ for any~$39\le\ell\le(f(\beta)m^*_1(f)-\delta)3n$ divisible by~$3$ or a blue tight cycle of length~$\ell$ for any~$39\le\ell\le\beta(f(\beta)m^*_1(f)-\delta)3n$ divisible by~$3$.
        \item \label{item:monox-cycle-embedding-xdiv-xdiv} $K$ contains a red tight cycle of length~$\ell$ for any~$39\le\ell\le(f(\beta)m^*_2(f)-\delta)3n$ or a blue tight cycle of length~$\ell$ for any~$39\le\ell\le\beta(f(\beta)m^*_2(f)-\delta)3n$.
        \item \label{item:monox-cycle-embedding-div-xdiv} $K$ contains a red tight cycle of length~$\ell$ for any~$39\le\ell\le(f(\beta)m^*_3(f)-\delta)3n$ divisible by~$3$ or a blue tight cycle of length~$\ell$ for any~$39\le\ell\le\beta(f(\beta)m^*_3(f)-\delta)3n$.
        \item \label{item:monox-cycle-embedding-xdiv-div} $K$ contains a red tight cycle of length~$\ell$ for any~$39\le\ell\le(f(\beta)m^*_4(f)-\delta)3n$ or a blue tight cycle of length~$\ell$ for any~$39\le\ell\le\beta(f(\beta)m^*_4(f)-\delta)3n$ divisible by~$3$.
    \end{enumerate}
\end{lemma}

We are now ready to prove Theorem~\ref{thm:main}\ref{item:thm-upper-bound}.

\begin{proof}[Proof of Theorem~\ref{thm:main}\ref{item:thm-upper-bound}]
    Let~$\eps,\delta,\eta>0$ and~$s,t\in\NN$ satisfy~$\frac{1}{s} \ll \eps \ll \delta \ll \eta$ and~$13\le t\le s$; set~$\beta:=\frac{t+1}{s+1}$. We distinguish four cases; the proof takes exactly the same form in each case and differs only in the parameters involved.

    \textbf{Case 1:}~$i=j=0$. Let~$K$ be a~$2$-edge-coloured~$3$-graph on~$N_1=3s+t+\eta s$ vertices. Note that by Lemma~\ref{lem:matching-div-div} we have~$m^*_1(\mu^*_1)=1$, so Lemma~\ref{lem:monox-cycle-embedding}\ref{item:monox-cycle-embedding-div-div} applied with~$\mu^*_1$ implies that~$K$ contains a red tight cycle of length~$\ell$ for any~$39\le\ell\le(\mu^*_1(\beta)-\delta)3N_1$ divisible by~$3$ or a blue tight cycle of length~$\ell$ for any~$39\le\ell\le\beta(\mu^*_1(\beta)-\delta)3N_1$ divisible by~$3$. Now we have~$(\mu^*_1(\beta)-\delta)3N_1 \ge 3(s+1)$, so~$K$ contains a red tight cycle of length~$3s$ or a blue tight cycle of length~$3t$ as required.

    \textbf{Case 2:}~$i,j\ne0$. Let~$K$ be a~$2$-edge-coloured~$3$-graph on~$N_2=\max(5s,3(s+t)) + \eta s$ vertices. Note that by Lemma~\ref{lem:matching-xdiv-xdiv} we have~$m^*_2(\mu^*_2)=1$, so Lemma~\ref{lem:monox-cycle-embedding}\ref{item:monox-cycle-embedding-xdiv-xdiv} applied with~$\mu^*_2$ implies that~$K$ contains a red tight cycle of length~$\ell$ for any~$39\le\ell\le(\mu^*_2(\beta)-\delta)3N_2$ or a blue tight cycle of length~$\ell$ for any~$39\le\ell\le\beta(\mu^*_2(\beta)-\delta)3N_2$. Now we have~$(\mu^*_2(\beta)-\delta)3N_2 \ge 3(s+1)$, so~$K$ contains a red tight cycle of length~$3s+i$ or a blue tight cycle of length~$3t+j$ as required.

    \textbf{Case 3:}~$i=0,\ j\ne0$. Let~$K$ be a~$2$-edge-coloured~$3$-graph on~$N_3=5s + \eta s$ vertices. Note that by Lemma~\ref{lem:matching-div-xdiv} we have~$m^*_3(\mu^*_3)=1$, so Lemma~\ref{lem:monox-cycle-embedding}\ref{item:monox-cycle-embedding-div-xdiv} applied with~$\mu^*_3$ implies that~$K$ contains a red tight cycle of length~$\ell$ for any~$39\le\ell\le(\mu^*_3(\beta)-\delta)3N_3$ divisible by~$3$ or a blue tight cycle of length~$\ell$ for any~$39\le\ell\le\beta(\mu^*_3(\beta)-\delta)3N_3$. Now we have~$(\mu^*_3(\beta)-\delta)3N_3 \ge 3(s+1)$, so~$K$ contains a red tight cycle of length~$3s$ or a blue tight cycle of length~$3t+j$ as required.

    \textbf{Case 4:}~$i\ne0,\ j=0$. Let~$K$ be a~$2$-edge-coloured~$3$-graph on~$N_4=\max(3s+t,5t) + \eta s$ vertices. Note that by Lemma~\ref{lem:matching-xdiv-div} we have~$m^*_4(\mu^*_4)=1$, so Lemma~\ref{lem:monox-cycle-embedding}\ref{item:monox-cycle-embedding-xdiv-div} applied with~$\mu^*_4$ implies that~$K$ contains a red tight cycle of length~$\ell$ for any~$39\le\ell\le(\mu^*_4(\beta)-\delta)3N_4$ or a blue tight cycle of length~$\ell$ for any~$39\le\ell\le\beta(\mu^*_4(\beta)-\delta)3N_4$ divisible by~$3$. Now we have~$(\mu^*_4(\beta)-\delta)3N_4 \ge 3(s+1)$, so~$K$ contains a red tight cycle of length~$3s+i$ or a blue tight cycle of length~$3t$ as required.
\end{proof}

\section{Proving Lemmas~\ref{lem:matching-xdiv-xdiv}--\ref{lem:matching-xdiv-div}} \label{sec:matching-xdivs}

In this section, we build up to the proofs of Lemmas~\ref{lem:matching-xdiv-xdiv}--\ref{lem:matching-xdiv-div}. In Section~\ref{ssec:setup} we introduce some structural analysis from~\cite{HaxellLuczakPengRodl} and encapsulate it in Setup~\ref{setup}. In Section~\ref{ssec:techniques-matching-colour-pivot} we develop several useful techniques and arguments which generally take a connected matching in one colour subject to some additional conditions and produce a connected matching in the other colour. Finally, we provide the proofs of Lemmas~\ref{lem:matching-xdiv-xdiv}--\ref{lem:matching-xdiv-div} in Section~\ref{ssec:proofs-lemmas}.

\subsection{The setup} \label{ssec:setup}

In this subsection, we introduce some structural analysis from~\cite{HaxellLuczakPengRodl}. We begin with two lemmas about the monochromatic components of~$2$-edge-coloured~$\eps$-complete~$3$-graphs. For completeness we provide proofs of these lemmas in Appendix~\ref{sec:lemma-proofs-matching-xdivs}.

\begin{lemma}[{\cite[Observation~8.2]{HaxellLuczakPengRodl}}] \label{lem:vertex-shadow-monoxcomp}
    Given a~$2$-edge-coloured~$\eps$-complete~$3$-graph~$H=H_{\red}\cup H_{\blue}$ on~$N$ vertices, for each vertex~$x\in V(H)$ there is a monochromatic component~$C_x$ such that
    \begin{equation} \label{eq:vertex-shadow-monoxcomp}
        |\{y\in V(H):xy\in \partial C_x\}|\ge(1-\eps)N\,.
    \end{equation}
\end{lemma}

\begin{lemma}[{\cite[Observation~8.4]{HaxellLuczakPengRodl}}] \label{lem:dominant-monoxcomp}
    Let~$H=H_{\red}\cup H_{\blue}$ be a~$2$-edge-coloured~$\eps$-complete~$3$-graph on~$N$ vertices. For each vertex~$x\in V(H)$ let~$C_x$ be a monochromatic component satisfying~\eqref{eq:vertex-shadow-monoxcomp}. Let~$V_{\red}$ denote the set of all vertices~$x\in V(H)$ such that~$C_x$ is red. If~$|V_{\red}|\ge6\eps N$, then there is a red component~$C_{\red}$ such that~$C_x=C_{\red}$ for all but at most~$2\eps N$ vertices~$x\in V_{\red}$.
\end{lemma}

Now we define the setup, which encapsulates the two lemmas above together with some additional definitions and properties; this serves as the setting of our proofs of Lemmas~\ref{lem:matching-div-div}--\ref{lem:matching-xdiv-div}. This encapsulation is convenient because we use this setup in many of the lemmas to come.

\begin{setup} \label{setup}
    When we say that we \emph{assume Setup~\ref{setup}} we mean that we have a~$2$-edge-coloured~$\eps$-complete~$3$-graph~$H=H_{\red}\cup H_{\blue}$ on~$N$ vertices such that for each vertex~$x\in V(H)$ there is a monochromatic component~$C_x$ satisfying~\eqref{eq:vertex-shadow-monoxcomp} and furthermore we have the following definitions. Set
    \[V_{\red}:=\{x\in V(H):C_x\textrm{ is red}\}\;\textrm{ and }\;V_{\blue}:=\{x\in V(H):C_x\textrm{ is blue}\}\,.\]
    If~$|V_{\red}|\ge6\eps N$, let the~\emph{scarlet} component~$S$ be the (unique) red component~$C_{\red}$ such that~$C_x=C_{\red}$ for all but at most~$2\eps N$ vertices~$x\in V_{\red}$ and set~$V'_{\red}=\{x\in V_{\red}:C_x=S\}$; if~$|V_{\red}|<6\eps N$, we say that~$S$ does not exist and set~$V'_{\red}=\nth$. If~$|V_{\blue}|\ge6\eps N$, let the~\emph{azure} component~$A$ be the (unique) blue component~$C_{\blue}$ such that~$C_x=C_{\blue}$ for all but at most~$2\eps N$ vertices~$x\in V_{\blue}$ and set~$V'_{\blue}=\{x\in V_{\blue}:C_x=A\}$; if~$|V_{\blue}|<6\eps N$, we say that~$A$ does not exist and set~$V'_{\blue}=\nth$. Set~$V':=V'_{\red}\cup V'_{\blue}$. For each~$x\in V'_{\red}$ set~$\partial S(x):=\{y\in V(H):xy\in\partial S\}$, and for each~$z\in V'_{\blue}$ set~$\partial A(z):=\{y\in V(H):zy\in\partial A\}$.
\end{setup}

The following lemma collects some useful properties about Setup~\ref{setup}.

\begin{lemma} \label{lem:setup-properties}
    In Setup~\ref{setup} with~$\eps<\frac{1}{12}$ the following hold.
    \begin{enumerate}[label=(\roman{*})]
        \item \label{item:complete-active} For all~$x\in V(H)$ we have~$|N_{\partial H}(x)|\ge(1-\eps)N+1$.
        \item \label{item:vertices-colour-robust} If~$|V_{\red}|\ge6\eps N$ then~$|V'_{\red}|\ge|V_{\red}|-2\eps N$; if~$|V_{\blue}|\ge6\eps N$ then~$|V'_{\blue}|\ge|V_{\blue}|-2\eps N$.
        \item \label{item:colour-bipartition-near-complete} $|V'|\ge N(1-8\eps)$.
        \item \label{item:dominant-monoxcomp-shadows-nbrhood-sizes} For all~$x\in V'_{\red}$ and~$z\in V'_{\blue}$ we have~$|\partial S(x)|,|\partial A(z)|\ge(1-\eps)N$.
    \end{enumerate}
\end{lemma}

\begin{proof}
    For~\ref{item:complete-active} let~$x\in V(H)$ and take~$y\in N_{\partial H}(x)\ne\nth$. Since~$\{y\}\cup N_H(x,y)\subseteq N_{\partial H}(x)$, the desired bound follows from~\eqref{eq:complete-codegree}. For~\ref{item:vertices-colour-robust} note that the desired outcome follows from Lemma~\ref{lem:dominant-monoxcomp} and the definitions of~$V'_{\red}$ and~$V'_{\blue}$. For~\ref{item:colour-bipartition-near-complete} note that~$\eps<\frac{1}{12}$ implies that we must have~$|V_{\red}|\ge6\eps N$ or~$|V_{\blue}|\ge6\eps N$; the desired outcome follows by combining this with~\ref{item:vertices-colour-robust}. For~\ref{item:dominant-monoxcomp-shadows-nbrhood-sizes} the desired outcome follows from Lemma~\ref{lem:vertex-shadow-monoxcomp} and the definitions of~$S$ and~$A$.
\end{proof}

\subsection{Constructing connected matchings} \label{ssec:techniques-matching-colour-pivot}

In this subsection, we collect several lemmas which encapsulate various techniques and arguments for constructing connected matchings. In particular, Lemmas~\ref{lem:matching-pivot-no-C4}--\ref{lem:not-scarlet-matching-force-blue-C4-matching} encapsulate various mechanisms which take a connected matching in a monochromatic component and produce a connected matching in a monochromatic component of the opposite colour; clearly, the colours red and blue may be freely interchanged. We begin with technical lemma which conveniently encapsulates a useful greedy-type argument.

\begin{lemma}[{\cite[Fact~A.2]{HaxellLuczakPengRodl2007}}] \label{lem:bipartite-high-degree-shared-edge}
    Let~$B$ and~$B'$ be two bipartite graphs on the same vertex set~$X\cup Y$. Suppose that there exist~$b,b'\in\NN$ such that
    \begin{enumerate} [label=(\roman{*})]
        \item for all~$x\in X$ we have~$\deg_{B}(x)\ge|Y|-b'$,
        \item for all~$y\in Y$ we have~$\deg_{B'}(y)\ge|X|-b$, and
        \item $(|X|-b)(|Y|-b')>bb'$.
    \end{enumerate}
    Then there exist vertices~$x\in X$ and~$y\in Y$ such that~$xy\in E(B)\cap E(B')$.
\end{lemma}

\begin{proof}
    Since the assumptions give~$|E(B)|+|E(B')|\ge|X|(|Y|-b')+(|X|-b)|Y|>|X||Y|$, the lemma follows.
\end{proof}

The following lemma encapsulates a technique which utilises the absence of copy of~$C^{(3)}_4$ in the monochromatic component containing the given connected matching and the presence of a disjoint set of vertices in the opposite colour; this technique is also utilised in~\cite{HaxellLuczakPengRodl}.

\begin{lemma} \label{lem:matching-pivot-no-C4}
    Assume Setup~\ref{setup} and suppose that~$C_{\red}$ is a red component such that~$C^{(3)}_4\nsubseteq C_{\red}$. Then for every connected matching~$M$ in~$C_{\red}$ and every subset~$W\subseteq V'_{\blue}\sm V(M)$ there is a connected matching~$M'$ of size~$\max(\min(|M|,|W|)-5\eps N,0)$ in~$A$ with an injective function~$\phi\colon M'\to M$ such that for each~$e\in M'$ we have~$|e\cap\phi(e)|=2$ and~$e\sm\phi(e)\subseteq W$.
\end{lemma}

\begin{proof}
    Fix a connected matching~$M$ in~$C_{\red}$ and a subset~$W\subseteq V'_{\blue}\sm V(M)$. If~$\min(|M|,|W|)\le5\eps N$, the empty matching suffices. Otherwise, we greedily construct~$M'$ and~$\phi\colon M'\to M$. Take~$M^*\subseteq M$ and~$W^*\subseteq W$ with~$|M^*|,|W^*|\ge 5\eps N$. Let~$B$ be the bipartite graph on~$M^*\cup W^*$ such that there is an edge between~$e=xyz\in M^*$ and~$w\in W^*$ if and only if~$xyw,xzw,yzw\in E(H)$; let~$B'$ be the bipartite graph on~$M^*\cup W^*$ such that there is an edge between~$e=xyz\in M^*$ and~$w\in W^*$ if and only if~$xw,yw,zw\in\partial A$. By~\eqref{eq:complete-codegree} we have~$\deg_{B}(e)\ge|W^*|-3\eps N$ for all~$e\in M^*$; by Lemma~\ref{lem:setup-properties}\ref{item:dominant-monoxcomp-shadows-nbrhood-sizes} we have~$\deg_{B'}(w)\ge|M^*|-\eps N$ for all~$w\in W^*$. Since~$(|M^*|-\eps N)(|W^*|-3\eps N)>3\eps^2N^2$, by Lemma~\ref{lem:bipartite-high-degree-shared-edge} there is an edge~$e=xyz\in M^*$ and a vertex~$w\in W^*$ such that~$xyw,xzw,yzw\in E(H)$ and~$xw,yw,zw\in\partial A$. Since~$C^{(3)}_4\nsubseteq C_{\red}$, at least one of the edges~$xyw,xzw,yzw$ must be blue and so in~$A$; pick one such edge and denote it by~$e'$. Since~$|e'\cap e|=2$ and~$e'\sm e\subseteq W$, we use~$e'$ to enlarge~$M'$ and extend~$\phi$ by setting~$\phi(e')=e$; this uses just one edge of~$M^*\subseteq M$ and one vertex of~$W^*\subseteq W$, so we may repeat this procedure until $M'$ has size~$\min(|M|,|W|)-5\eps N$.
\end{proof}

The following lemma explores a different technique which exploits a maximality condition on the number of vertices of the same colour as the given matching and utilises the edges which have many opposite-colour vertices; this technique is an innovation which allows us to attain the~$5s$-type quantities observed in our lower bound constructions instead of the the~$6s$-type quantities seen in~\cite{HaxellLuczakPengRodl}.

\begin{lemma} \label{lem:matching-pivot-max-red-vertices}
    Assume Setup~\ref{setup}. Let~$C_{\red}$ be a red component, $U\subseteq V(H)$ and~$s\in\NN$. Let~$\cM$ be the set of all matchings of size~$s$ in~$C_{\red}$ whose vertex set is contained in~$U$. Suppose that~$M\in\cM$ has the maximum number of vertices in~$V'_{\red}$ among all elements of~$\cM$. Let~$M^*:=\{e\in M:|e\cap V'_{\blue}|\ge2\}$. Then for every subset~$W\subseteq (U\cap V'_{\red})\sm V(M)$ there is a matching~$M'$ of size~$\min(|M^*|,\max(|W|-2\eps N,0))$ in~$A$ with an injective function~$\phi\colon M'\to M^*$ such that for each~$e\in M'$ we have~$|e\cap\phi(e)|=2$, $e\sm\phi(e)\subseteq W$ and~$\phi(e)\sm e\subseteq V'_{\blue}$.
\end{lemma}

\begin{proof}
    For each edge~$e\in M^*$ fix a pair~$x_e,y_e\in e$ of distinct vertices such that~$e\sm\{y_e\}\subseteq V'_{\blue}$. Let~$B$ be the bipartite graph on~$M^*\cup W$ such that there is an edge between~$e\in M^*$ and~$w\in W$ if and only if~$x_ey_ew\in E(H)$ and~$x_ew\in\partial A$. By~\eqref{eq:complete-codegree} and Lemma~\ref{lem:setup-properties}\ref{item:dominant-monoxcomp-shadows-nbrhood-sizes} we have~$\deg_B(e)\ge|W|-2\eps N$ for all~$e\in M^*$, so $B$ contains a matching~$M_B$ of size~$\min(|M^*|,\max(|W|-2\eps N,0))$. Observe that for each~$ew\in M_B$ the edge~$x_ey_ew\in E(H)$ is a blue edge in~$A$. Indeed, if~$x_ey_ew$ were red then~$(M\sm\{e\})\cup\{x_ey_ew\}\in\cM$ would have more vertices in~$V'_{\red}$ than~$M$, contradicting the maximality of~$M$. Set~$M':=\{x_ey_ew:ew\in M_B\}$ and define~$\phi\colon M'\to M^*$ as follows. For each~$f\in M'$ set~$\phi(f)$ to be the unique~$e\in M^*$ such that~$x_e,y_e\in f$.
\end{proof}

The following lemma encapsulates a technique which uses the assumption that the given matching is not in a special component to produce a connected matching in a monochromatic component of the opposite colour containing a copy of~$C^{(3)}_4$. While a very specific version of this argument with~$V(M)\subseteq V'_{\red}$ features in~\cite{HaxellLuczakPengRodl}, our key innovation lies in handling the non-$V'_{\red}$ vertices via the partition~$\{P'_0,P'_1,P'_2,P'_3\}$ below.

\begin{lemma} \label{lem:not-scarlet-matching-force-blue-C4-matching}
    Assume Setup~\ref{setup}. Suppose that there is a connected matching~$M$ in a red component~$C_{\red}\ne S$. For~$i\in\{0,1,2,3\}$ set~$P'_i:=\{e\in M:|e\cap V'_{\red}|=i\}$. Suppose further that
    \[q := |M|(1-12\eps^{1/2})-\max(12\eps^{1/2}|M|+|P'_0|-|P'_2|-|P'_3|,0)-2\eps^{1/2}N-2 >0\]
    is positive. Then there is a blue component~$C_{\blue}$ containing a copy of~$C^{(3)}_4$ and a connected matching of size at least~$q$.
\end{lemma}

\begin{proof}
    We say that an ordered pair~$(e,f)\in M\times M$ of distinct edges is \emph{good} if for each~$x\in e\cap V'_{\red}$ and each~$z\in f$ we have~$xz\in\partial S$ and for each pair~$xy\subseteq e$ and each~$z\in f$ we have~$xyz\in E(H)$; call an unordered pair~$\{e,f\}\subseteq M$ of distinct edges \emph{compatible} if both~$(e,f)$ and~$(f,e)$ are good. Note that every compatible pair~$\{e,f\}$ satisfies the following.
    \begin{enumerate}[label=(PL\arabic{*})]
        \item \label{item:cross-complete-S-shadow} $H[e\cup f]\cong K^{(3)}_6$ and for all~$x\in e$ and~$z\in f$ with~$\{x,z\}\cap V'_{\red}\ne\nth$ we have~$xz\in\partial S$.
    \end{enumerate}
    Let~$M^*\subseteq M$ be the collection of edges~$e\in M$ such that there are at most~$\eps^{1/2}N$ edges~$f\in M$ for which the pair~$\{e,f\}$ is not compatible. By~\eqref{eq:complete-codegree} and Lemma~\ref{lem:setup-properties}\ref{item:dominant-monoxcomp-shadows-nbrhood-sizes}, for each edge~$e\in M$ there are at most~$6\eps N$ edges~$f\in M$ such that~$(e,f)$ is not good. Hence, there are at most~$6\eps N|M|$ incompatible pairs and so
    \begin{equation} \label{eq:matching-compatible}
        |M^*|\ge(1-12\eps^{1/2})|M|\,.
    \end{equation}
    Set~$P_0:=P'_0\cap M^*$, $P_1:=P'_1\cap M^*$ and~$P_2:=(P'_2\cup P'_3)\cap M^*$. Since~$q>0$, we have~$P_1\cup P_2\ne\nth$. Observe that by~\eqref{eq:matching-compatible} we have
    \begin{equation} \label{eq:size-ineq}
        |P_0|-|P_2|\le|P'_0|-|P'_2|-|P'_3|+12\eps^{1/2}|M|\,.
    \end{equation}
    
    We shall construct a collection~$\cC=\{e_1,\dots,e_q,e_{q+1}\}\subseteq M^*\subseteq M$ of edges such that for all~$\ell\in[q]$ the pair~$\{e_{\ell},e_{\ell+1}\}$ is compatible and satisfies~$|(e_{\ell}\cup e_{\ell+1})\cap V'_{\red}|\ge2$; the exact procedure depends on whether~$P_2\ne\nth$. If~$P_2=\nth$, start with an arbitrary~$e_1\in P_1\ne\nth$ and when we have chosen~$\cC_{\ell}=\{e_1,\dots,e_{\ell}\}$, arbitrarily pick~$e_{\ell+1}\in P_1\sm\cC_\ell$ such that~$\{e_{\ell},e_{\ell+1}\}$ is compatible. This procedure stops only when at most~$\eps^{1/2}N$ edges of~$P_1$ remain unused, so~$\cC$ contains at least~$|M^*|-|P_0|-\eps^{1/2}N\ge q+1$ edges. If~$P_2\ne\nth$, start with an arbitrary~$e_1\in P_2$ and proceed in two stages. Suppose that we have chosen~$\cC_{\ell}=\{e_1,\dots,e_{\ell}\}$ with~$e_{\ell}\in P_{i_{\ell}}$; an edge~$e\in M^*$ is \emph{valid} if the pair~$e_{\ell}$ and~$e$ is compatible. If we are in stage~$1$, arbitrarily pick a valid~$e_{\ell+1}\in P_{2-i_{\ell}}\sm\cC_\ell$ if possible. If not, proceed to stage~$2$; roll back to~$e_{\ell-1}$ if~$i_{\ell}=0$. If we are in stage~$2$, arbitrarily pick a valid~$e_{\ell+1}\in (P_1\cup P_2)\sm\cC_\ell$ if possible. Observe that there are at most~$\eps^{1/2}N$ unused edges in~$P_1\cup P_2$ when the procedure stops; since the edges chosen in stage~$1$ alternate between~$P_2$ and~$P_0$, there are at most~$\max(|P_0|-|P_2|,0)+\eps^{1/2}N+1$ unused edges in~$P_0$ when the procedure stops. Hence, by~\eqref{eq:matching-compatible} and~\eqref{eq:size-ineq} $\cC$ contains at least
    \[ |M^*| - \max(|P_0|-|P_2|,0)-2\eps^{1/2}N-1 \ge q+1 \]
    edges as required.

    Let~$\ell\in[q]$ and let~$F_{\ell}:=\{f\in H[e_{\ell}\cup e_{\ell+1}]:f\cap e_{\ell},f\cap e_{\ell+1},f\cap V'_{\red}\ne\nth\}$ be the set of cross-edges on~$e_{\ell}\cup e_{\ell+1}$ which meet~$V'_{\red}$. Since edges in~$F_{\ell}$ all contain a pair in~$\partial S$ by~\ref{item:cross-complete-S-shadow} for~$\{e_{\ell},e_{\ell+1}\}$ and a pair in~$\binom{e_{\ell}}{2}\cup\binom{e_{\ell+1}}{2}\subseteq\partial C\ne\partial S$, it follows that they must all be blue and belong to a single blue component. Indeed, a red edge in~$F_{\ell}$ is in both~$C_{\red}$ and~$S$, which gives a contradiction. Since for each~$\ell\in[q-1]$ there is a pair~$xy\in\binom{e_{\ell+1}}{2}$ such that each of~$F_{\ell}$ and~$F_{\ell+1}$ contains an edge which contains~$xy$, it follows that the edges in~$F_{\ell}\cup F_{\ell+1}$ belong to a single blue component. Hence, the edges in~$F:=\bigcup_{\ell\in[s]}F_i$ belong to a single blue component~$C_{\blue}$. Pick~$S\subseteq e_1\cup e_2$ with~$|S\cap e_1|=|S\cap e_2|=|S\cap(e_1\cup e_2)\cap V'_{\red}|=2$; this is possible because~$|(e_1\cup e_2)\cap V'_{\red}|\ge2$. Since~$H[S]$ is a copy of~$C^{(3)}_4$ by~\ref{item:cross-complete-S-shadow} and~$E(H[S])\subseteq F_1$, this yields a copy of~$C^{(3)}_4$ in~$C_{\blue}$. Now for each~$1\le i\le \left\lceil\frac{q}{2}\right\rceil$ fix distinct~$u_i,v_i\in(e_{2i-1}\cup e_{2i})\cap V'_{\red}$ and pick a pair~$f_{2i-1},f_{2i}\in F_{2i-1}$ of disjoint edges such that~$u_i\in f_{2i-1}$ and~$v_i\in f_{2i}$; this yields a connected matching~$\{f_1,\dots,f_q\}$ of size at least~$q$ in~$C_{\blue}$.
\end{proof}

The following lemma encapsulates a key argument in the proofs of Lemmas~\ref{lem:matching-xdiv-xdiv}--\ref{lem:matching-xdiv-div}. Broadly speaking, this builds on the simple greedy argument at the end of Subcase 2A in~\cite[Lemma 2.2]{HaxellLuczakPengRodl2007} and utilises an innovative technique to reduce the number of vertices required from the~$6s$-type quantities seen in~\cite{HaxellLuczakPengRodl} to the~$5s$-type quantities observed in our lower bound constructions.

\begin{lemma} \label{lem:build-blue-matching-C4-small-red-set-blue-C4-dominant-pair}
    Assume Setup~\ref{setup} with~$1/N \ll \eps \ll 1$. Suppose that
    \begin{enumerate}[label=(BR\arabic{*})]
        \item \label{item:few-blue-vertices} $\left|V'_{\blue}\right| < T := \frac{N}{5}(1-26\eps^{1/12})$,
        \item \label{item:red-scarlet-matching} there is a connected matching~$M_S$ of size~$|V'_{\blue}|$ in~$S$ with~$V(M_S)\subseteq V'_{\red}$,
        \item \label{item:scarlet-no-C4} $C^{(3)}_4\nsubseteq S$, and
        \item \label{item:small-red-set-blue-C4-dominant-pair} there exist two blue components~$C^{(1)}_{\blue}$ and~$C^{(2)}_{\blue}$ such that~$C^{(1)}_{\blue}=A$ whenever~$A$ exists and every set~$X\subseteq V'_{\red}$ with~$|X|\ge25\eps^{1/12}N$ contains a copy of~$C^{(3)}_4$ in~$C^{(1)}_{\blue}$ or~$C^{(2)}_{\blue}$.
    \end{enumerate}
    Then there is a blue component~$C_{\blue}$ containing a copy of~$C^{(3)}_4$ and a connected matching of size~$T$.
\end{lemma}

\begin{proof}
    Fix two blue components~$C^{(1)}_{\blue}$ and~$C^{(2)}_{\blue}$ returned by~\ref{item:small-red-set-blue-C4-dominant-pair}. We iteratively construct two matchings~$M_1$ and~$M_2$ as follows. If~$V'_{\blue}=\nth$, we start with~$M_1=M_2=\nth$. If~$V'_{\blue}\ne\nth$, then~$A$ exists and Lemma~\ref{lem:matching-pivot-no-C4} applied with~\ref{item:scarlet-no-C4}, $C_{\red}=S$, the connected matching~$M_S$ given by~\ref{item:red-scarlet-matching} and~$W=V'_{\blue}$ returns a connected matching~$M'$ of size at least~$|V'_{\blue}|-5\eps N$ in~$A$ with an injective function~$\phi\colon M'\to M_S$ such that for each~$e\in M'$ we have~$|e\cap\phi(e)|=2$ and~$e\sm\phi(e)\subseteq V'_{\blue}$; in this case, we start with~$M_1=M'$ and~$M_2=\nth$. Throughout the iterative procedure, we will maintain the following.
    \begin{enumerate}[label=(M\arabic{*})]
        \item \label{item:matchings-intersection-pairing} There exist~$M^*_1\subseteq M_1$ and~$M^*_2\subseteq M_2$ with~$|M^*_1|=|M^*_2|\ge\min(|M_1|,|M_2|)-5\eps N$ and a bijection~$\beta\colon M^*_1\to M^*_2$ such that for all~$e_1\in M_1$ and~$e_2\in M_2$ we have~$e_1\cap e_2\ne\nth$ if and only if~$e_1\in M^*_1$, $e_2\in M^*_2$ and~$\beta(e_1)=e_2$; furthermore, $M^*_1$, $M^*_2$ and~$\beta$ are unique.
        \item \label{item:blue-fully-utilised} At most~$5\eps N$ vertices of~$V'_{\blue}$ are not in $V(M_1)\cup V(M_2)$.
        \item \label{item:matching-edges-blue-vertices-restricted} If~$e\in M_1\cup M_2$ satisfies~$e\cap V'_{\blue}\ne\nth$, then~$e\notin M_2\sm M^*_2$ and~$|e\cap V'_{\blue}|=1$.
    \end{enumerate}
    Observe that these are initially all true with~$M^*_1=M^*_2=\nth$.

    If the set~$Y:=V'_{\red}\sm(V(M_1)\cup V(M_2))$ satisfies~$|Y|\ge25\eps^{1/12}N$, we may fix~$Y^*\subseteq Y$, $M^*_1\subseteq M_1$, $M^*_2\subseteq M_2$ and~$\beta\colon M^*_1\to M^*_2$ satisfying~$|Y^*|=25\eps^{1/12}N$ and~\ref{item:matchings-intersection-pairing}. By~\ref{item:small-red-set-blue-C4-dominant-pair}~$Y^*$ contains a copy~$C$ of~$C^{(3)}_4$ in~$C^{(1)}_{\blue}$ or~$C^{(2)}_{\blue}$; fix an edge~$e$ of~$C$. Since the arguments are very similar, we shall assume that~$e\in C^{(1)}_{\blue}$. If there exist~$\{x_e,y_e\}\subseteq e$ and~$\{x_f,y_f\}\subseteq f\in M_2\sm M^*_2$ such that~$Q:=H[\{x_e,y_e,x_f,y_f\}]$ is a copy of~$C^{(3)}_4$, the bipartite induced subgraph~$\partial S[\{x_e,y_e\},\{x_f,y_f\}]$ is a copy of~$K_{2,2}$ and~$(e\sm\{x_e,y_e\})\cap V'_{\blue}=\nth$, then by~\ref{item:scarlet-no-C4} there is a blue edge~$g\in E(Q)$; if~$|g\cap e|=2$, then add~$g$ to~$M_1$ and~$M^*_1$, add~$f$ to~$M^*_2$ and extend~$\beta$ by setting~$\beta(g)=f$; if~$|g\cap f|=2$, then add~$e$ to~$M_1$ and~$M^*_1$, swap~$f$ for~$g$ in~$M_2$, add~$g$ to~$M^*_2$ and extend~$\beta$ by setting~$\beta(e)=g$. Otherwise, add~$e$ to~$M_1$ and keep~$M^*_1$, $M_2$, $M^*_2$ and~$\beta$ unchanged.
    
    Observe that each iteration maintains the fact that~$M_1$ and~$M_2$ are matchings and increases the quantity~$|(V(M_1)\cup V(M_2)|$ by at least two. To see that we maintain~\ref{item:matchings-intersection-pairing}--\ref{item:matching-edges-blue-vertices-restricted}, it suffices to verify that the former case always applies when~$|M^*_1|=|M^*_2|=\min(|M_1|,|M_2|)-5\eps N$; the required outcomes then follow from the definitions and procedure. To this end, set~$X_1:=(M_1\sm M^*_1)\cup\{e\}$ and~$X_2:=M_2\sm M^*_2$; note that~$|X_1|,|X_2|\ge5\eps N$. For each~$e\in X_1\cup X_2$ fix an arbitrary pair~$\{x_e,y_e\}\subseteq e$ such that~$(e\sm\{x_e,y_e\})\cap V'_{\blue}=\nth$. Let~$B$ be the bipartite graph on~$X_1\cup X_2$ such that there is an edge between~$e\in X_1$ and~$f\in X_2$ if and only if~$x_ey_ex_f,x_ey_ey_f\in E(H)$; let~$B'$ be the bipartite graph on~$X_1\cup X_2$ such that there is an edge between~$e\in X_1$ and~$f\in X_2$ if and only if~$x_ex_fy_f,y_ex_fy_f\in E(H)$ and~$x_ex_f,x_ey_f,y_ex_f,y_ey_f\in\partial S$. By~\eqref{eq:complete-codegree} we have~$\deg_{B}(e)\ge|X_2|-\eps N$ for all~$e\in X_1$; by~\eqref{eq:complete-codegree} and Lemma~\ref{lem:setup-properties}\ref{item:dominant-monoxcomp-shadows-nbrhood-sizes} we have~$\deg_{B'}(f)\ge|X_1|-3\eps N$ for all~$f\in X_2$. Since~$(|X_1|-3\eps N)(|X_2|-\eps N)>3\eps^2N^2$, by Lemma~\ref{lem:bipartite-high-degree-shared-edge} there exists~$f\in X_2$ such that~$\{x_e,y_e\}\subseteq e$ and~$\{x_f,y_f\}\subseteq f$ witnesses the former case.

    Since~\ref{item:matchings-intersection-pairing}--\ref{item:matching-edges-blue-vertices-restricted} and the applicability of~\ref{item:small-red-set-blue-C4-dominant-pair} means that we may repeat this procedure until~$M_1$ and~$M_2$ cover all but at most~$25\eps^{1/12}N + 5\eps N$ vertices of~$V'$ and $V(M_1)\cap V(M_2)$ contains at least~$\min(|M_1|,|M_2|)-5\eps N$ vertices, by Lemma~\ref{lem:setup-properties}\ref{item:colour-bipartition-near-complete} we have
    \[3(|M_1|+|M_2|) - \min(|M_1|,|M_2|) + 5\eps N \ge |V(M_1)\cup V(M_2)| \ge N(1 - 25\eps^{1/12} - 13\eps)\,,\]
    which implies that for some~$i\in[2]$ the matching~$M_i$ has at least~$T$ edges; note that by construction and size, at least one edge of~$M_i$ is derived from a copy of~$C^{(3)}_4$ in~$C^{(i)}_{\blue}$.
\end{proof}

The following lemma encapsulates an argument which we use repeatedly in the proofs of Lemmas~\ref{lem:matching-xdiv-xdiv}--\ref{lem:matching-xdiv-div}. Broadly speaking, this corresponds to Case 2 in~\cite[Lemma 2.2]{HaxellLuczakPengRodl2007} and replaces the simple greedy argument at the end of Subcase 2A with Lemma~\ref{lem:build-blue-matching-C4-small-red-set-blue-C4-dominant-pair}.

\begin{lemma} \label{lem:build-blue-C4-matching-scarlet-C5}
    Assume Setup~\ref{setup} with~$1/N \ll \eps \ll 1$ and suppose that~\ref{item:few-blue-vertices}--\ref{item:scarlet-no-C4} hold. Then there is a blue component~$C_{\blue}$ containing a copy of~$C^{(3)}_4$ and a connected matching of size~$T$, or $S$ contains a copy of $C^{(3)}_5$.
\end{lemma}

Lemma~\ref{lem:build-blue-C4-matching-scarlet-C5} follows immediately from the combination of Lemma~\ref{lem:build-blue-matching-C4-small-red-set-blue-C4-dominant-pair} and the following two lemmas. The arguments which yield these lemmas are available in~\cite{HaxellLuczakPengRodl2007}, and for completeness we provide proofs of these lemmas in Appendix~\ref{sec:lemma-proofs-matching-xdivs}.

\begin{lemma} \label{lem:dominant-pair-small-red-set-blue-C4}
    Assume Setup~\ref{setup} with~$1/N \ll \eps \ll 1$. Suppose that~\ref{item:few-blue-vertices}--\ref{item:scarlet-no-C4} hold and that there exist two blue components~$C^{(1)}_{\blue}$ and~$C^{(2)}_{\blue}$ such that
    \begin{equation} \label{eq:blue-monoxcomps-dominant-pair}
        \left|\binom{V'_{\red}}{2}\sm\left(\partial C^{(1)}_{\blue}\cup\partial C^{(2)}_{\blue}\right)\right| < 2\eps^{1/3}N^2\,.
    \end{equation}
    Then~\ref{item:small-red-set-blue-C4-dominant-pair} holds.
\end{lemma}

\begin{lemma} \label{lem:no-dominant-pair-scarlet-C5}
    Assume Setup~\ref{setup} with~$1/N \ll \eps \ll 1$. Suppose that~\ref{item:few-blue-vertices}--\ref{item:scarlet-no-C4} hold and that for any two blue components~$C^{(1)}_{\blue}$ and~$C^{(2)}_{\blue}$ the inequality~\eqref{eq:blue-monoxcomps-dominant-pair} does not hold. Then~$S$ contains a copy of $C^{(3)}_5$.
\end{lemma}

\begin{proof}[Proof of Lemma~\ref{lem:build-blue-C4-matching-scarlet-C5}]
    If there exist two blue components~$C^{(1)}_{\blue}$ and~$C^{(2)}_{\blue}$ such that~\eqref{eq:blue-monoxcomps-dominant-pair} holds, then~\ref{item:small-red-set-blue-C4-dominant-pair} holds by Lemma~\ref{lem:dominant-pair-small-red-set-blue-C4}, so by Lemma~\ref{lem:build-blue-matching-C4-small-red-set-blue-C4-dominant-pair} there is a blue component~$C_{\blue}$ containing a copy of~$C^{(3)}_4$ and a connected matching of size~$T$. Otherwise, \eqref{eq:blue-monoxcomps-dominant-pair} does not hold for any two blue components~$C^{(1)}_{\blue}$ and~$C^{(2)}_{\blue}$. In this case, by Lemma~\ref{lem:no-dominant-pair-scarlet-C5}~$S$ contains a copy of~$C^{(3)}_5$.
\end{proof}

\subsection{Proofs of Lemmas~\ref{lem:matching-xdiv-xdiv}--\ref{lem:matching-xdiv-div}} \label{ssec:proofs-lemmas}

\begin{proof}[Proof of Lemma~\ref{lem:matching-xdiv-xdiv}]
    By Lemmas~\ref{lem:vertex-shadow-monoxcomp} and~\ref{lem:dominant-monoxcomp}, from now onwards we may assume Setup~\ref{setup} with~$N=\max(5s,3(s+t))+\eta s$. Suppose that the~$2$-edge-colouring~$H=H_{\red}\cup H_{\blue}$ is such that
    \begin{equation} \label{eq:xdiv-xdiv-contra}
        \begin{split}
            &\textrm{no red component contains }\;M^{(3)}_s\;\textrm{ and either }\;C^{(3)}_4\;\textrm{ or }\;C^{(3)}_5\;\textrm{, and } \\
            &\textrm{no blue component contains }\;M^{(3)}_t\;\textrm{ and either }\;C^{(3)}_4\;\textrm{ or }\;C^{(3)}_5\,.
        \end{split}
    \end{equation}
    We distinguish two cases and obtain a contradiction to~\eqref{eq:xdiv-xdiv-contra} in each of them. Set~$\delta:=\frac{\eta}{20}$ and~$s':=(1+\delta)s$.

    \textbf{Case 1:}~$|V'_{\red}|,|V'_{\blue}|\ge s'$. We shall assume that~$|V'_{\red}|\ge|V'_{\blue}|$; the proof for~$|V'_{\red}|\le|V'_{\blue}|$ is identical with red and blue interchanged. Partition~$V'=U\cup U'$ with~$U'\subseteq V'_{\blue}$, $|U'|=s'$ and~$U=V'\sm U'$. We prove the following claim.

    \begin{claim} \label{claim:without-reservoir-majority-dominant-monoxcomp-matching}
        There is a matching of size~$s'$ in~$S$ whose vertex set is contained in~$U$.
    \end{claim}

    \begin{claimproof}
        Suppose that there is no matching of size~$s'$ in~$S$ whose vertex set is contained in~$U$, and suppose further that there is a matching~$M$ of size~$s'+\delta s$ in a blue component~$C$ of~$H$ whose vertex set is contained in~$U$; without loss of generality, let~$M$ have the maximum number of vertices in~$V'_{\blue}$ among all matchings of size~$s'+\delta s$ in~$C$ whose vertex set is contained in~$U$. Set~$W:=(U\cap V'_{\blue})\sm V(M)$, $X:=V'_{\red}\sm V(M)$ and~$M^*:=\{e\in M:|e\cap V'_{\red}|\ge2\}$; note that
        \begin{equation} \label{eq:M^*-size}
            |M^*|\ge\min\left(|M|,\max\left(\frac{|V'_{\red}|-|M|-|X|}{2},0\right)\right)\,.
        \end{equation}
        Apply Lemma~\ref{lem:matching-pivot-max-red-vertices} to~$U$, $M$ and~$W$ with the colours red and blue interchanged to obtain a matching~$M'$ of size~$\min(|M^*|,\max(|W|-2\eps N,0))$ in~$S$ with a function~$\phi'\colon M'\to M^*$. Set~$M_1:=M\sm\phi'(M')$ and~$X_1:=V'_{\red}\sm(V(M')\cup V(M_1))$; note that~$|X_1|=|X|+|M'|$. Now apply Lemma~\ref{lem:matching-pivot-no-C4} to~$M_1$ and~$X_1$ with the colours red and blue interchanged to obtain a matching~$M''$ of size~$\min(|M_1|,|X_1|) - 5\eps N$ in~$S$ with a function~$\phi''\colon M''\to M_1$. The properties of~$\phi'$ and~$\phi''$ imply that~$M'\cup M''$ is a connected matching in~$S$ of size
        \begin{align*}
            |M'| + \min(|M_1|,|X_1|) - 5\eps N & = \min(|M|,2|M'|+|X|) - 5\eps N \\
            & \ge \min(|M|,|V'_{\red}|-|M|,|U\sm V(M)|-2\eps N) - 5\eps N \ge s'\,,
        \end{align*}
        where we apply~\eqref{eq:M^*-size}, the fact that~$|V'_{\red}|\ge|V'_{\blue}|$ and Lemma~\ref{lem:setup-properties}\ref{item:colour-bipartition-near-complete}. But this is a contradiction, so no blue component of~$H$ contains a matching of size~$s'+\delta s$ whose vertex set is contained in~$U$.

        By applying Lemma~\ref{lem:matching-div-div} to the induced colouring~$H[U]=H_{\red}[U]\cup H_{\blue}[U]$, we obtain a matching~$M'$ of size~$\frac{|V'|-2s'-\delta s}{3+\delta}\ge s\left(1+\frac{\eta}{4}\right)$ in some red component~$C$ of~$H$ with~$V(M')\subseteq U$. For~$i\in[3]\cup\{0\}$ set~$P'_i:=\{e\in M':|e\cap V'_{\red}|=i\}$. Since~$|V'_{\red}|\ge|V'_{\blue}|$ and~$U'\subseteq V'_{\blue}$, we obtain~$|V'|/2-|U'| \ge |V'_{\blue}\cap V(M')| \ge 3|P'_0| + 2|P'_1| + |P'_2| \ge 3|M'|/2 + |P'_0|-|P'_2|-|P'_3|$, which yields
        \begin{equation} \label{eq:bound-red-imbalance}
            |P'_0|-|P'_2|-|P'_3|\le|V'|/2-3|M'|/2-|U'| = \delta|M'|/2 + \delta s/2\,.
        \end{equation}
        Suppose that~$C\ne S$. Then Lemma~\ref{lem:not-scarlet-matching-force-blue-C4-matching} applied to~$M'$ with~\eqref{eq:bound-red-imbalance} returns a blue component containing a copy of~$C^{(3)}_4$ and a matching of size~$|M'|\left(1-24\eps^{1/2}-\delta/2\right) - \delta s/2 - 2\eps^{1/2}N - 2 \ge s$, which contradicts~\eqref{eq:xdiv-xdiv-contra}.
    \end{claimproof}

    Claim~\ref{claim:without-reservoir-majority-dominant-monoxcomp-matching} returns a matching~$M$ of size~$s'\ge s$ in~$S$ with~$V(M)\subseteq U$ and Lemma~\ref{lem:matching-pivot-no-C4} applied to~$M$ and~$U'$ returns a matching of size at least~$s'-5\eps N \ge s$ in~$A$, so by~\eqref{eq:xdiv-xdiv-contra} we have
    \begin{equation} \label{eq:dominant-monoxcomps-no-C4}
        C^{(3)}_4\nsubseteq S\;\textrm{ and }\;C^{(3)}_4\nsubseteq A\,.
    \end{equation}
    We use~\eqref{eq:dominant-monoxcomps-no-C4} and the fact that~$r(C^{(3)}_4)=13$ (see~\cite{McKayRadziszowski}) to prove the following claim.

    \begin{claim} \label{claim:all-blue-monoxcomp}
        Every edge~$e\in E(H)\cap\binom{V'_{\red}}{3}$ is blue; consequently~$E(H)\cap\binom{V'_{\red}}{3}\subseteq C'_{\blue}$ for some blue component~$C'_{\blue}$.
    \end{claim}

    \begin{claimproof}
        We begin with the following subclaim.

        \begin{claim} \label{claim:pair-few-blue-cross-edges}
            For each~$xy\in\partial H\cap \binom{V'_{\red}}{2}$ the set~$Z_{xy}:=\{z\in V'_{\blue}:xyz\in E(H_{\blue})\}$ of vertices satisfies~$|Z_{xy}|\le80\eps N$.
        \end{claim}

        \begin{claimproof}
            Suppose there exists~$xy\in\partial H\cap \binom{V'_{\red}}{2}$ such that~$|Z_{xy}|>80\eps N$. We shall pick a set of thirteen distinct vertices~$Z=\{z_1,\dots,z_{13}\}\subseteq Z_{xy}$ such that
            \begin{enumerate}[label=(\roman{*})]
                \item \label{item:complete-edges} for all~$1\le i<j<k\le 13$ we have~$z_iz_jz_k\in E(H)$ and
                \item \label{item:both-dominant-monoxcomps} for all~$1\le i<j\le 13$ we have~$z_iz_j\in\partial S\cap\partial A$.
            \end{enumerate}
            Pick~$z_1\in Z$ so that~$xz_1,yz_1\in\partial S$; by Lemma~\ref{lem:setup-properties}\ref{item:dominant-monoxcomp-shadows-nbrhood-sizes}, at most~$2\eps N$ vertices of~$Z$ are not valid. For each~$i=2,\dots,13$ in turn pick~$z_i\in Z$ so that~$xz_i,yz_i\in\partial S$, $z_iz_j\in\partial A$ for~$j=1,\dots,i-1$ and~$z_iz_jz_k\in E(H)$ for all~$1\le j<k<i$; since by~\eqref{eq:complete-codegree} and Lemma~\ref{lem:setup-properties}\ref{item:dominant-monoxcomp-shadows-nbrhood-sizes} at most~$2\eps N + (i-1)\eps N + \binom{i-1}{2}\eps N \le 80\eps N<|Z|$ vertices are not valid, we may pick~$z_i$. The vertices~$z_1,\dots,z_{13}$ clearly satisfy~\ref{item:complete-edges}. Furthermore, for~$1\le i<j\le13$ we have~$z_iz_j\in\partial A$ and the edges~$xyz_1,\dots,xyz_{13}$ are blue. Take~$1\le i<j\le13$. Since~$C^{(3)}_4\nsubseteq A$ by~\eqref{eq:dominant-monoxcomps-no-C4} and~$xz_i,xz_j,yz_i,yz_j\in\partial S$, one of the edges~$xz_iz_j$ and~$yz_iz_j$ must be a red edge in~$S$. Hence, we have~$z_iz_j\in\partial S$ and~\ref{item:both-dominant-monoxcomps} holds.

            Since~$r(C^{(3)}_4)=13$, there is a monochromatic copy~$C$ of~$C^{(3)}_4$ in~$H[Z]$. Now if~$C$ were red then~$C\subseteq S$ and if~$C$ were blue then~$C\subseteq A$; both outcomes contradict~\eqref{eq:dominant-monoxcomps-no-C4}.
        \end{claimproof}

        Fix~$xyz\in E^*:=E(H)\cap\binom{V'_{\red}}{3}$. By Claim~\ref{claim:pair-few-blue-cross-edges} there are at most~$240\eps N$ vertices~$a\in V'_{\blue}$ such that at least one of the edges~$xya,xza,yza$ is blue; by~\eqref{eq:complete-codegree} there are at most~$3\eps N$ vertices~$a\in V'_{\blue}$ such that at least one of the triples~$xya,xza,yza$ is not an edge of~$H$; by Lemma~\ref{lem:setup-properties}\ref{item:dominant-monoxcomp-shadows-nbrhood-sizes} there are at most~$3\eps N$ vertices~$a\in V'_{\blue}$ such that at least one of the pairs~$xa,ya,za$ is not in~$\partial S$. Since~$|V'_{\blue}|\ge\left(1+\frac{\eta}{20}\right)s > 246\eps N$, we may pick a vertex~$a\in V'_{\blue}$ such that~$xya,xza,yza\in S$. Since~$C^{(3)}_4\nsubseteq S$ by~\eqref{eq:dominant-monoxcomps-no-C4}, the edge~$xyz$ must be blue.

        To see that the edges in~$E^*$ are all in the same blue component, first note that any pair of blue edges in~$E^*$ which share two vertices are clearly in the same blue component. Now consider two edges~$z_1z_2z_3$ and~$z_3z_4z_5$ in~$E^*$ which share one vertex. Since~$|V'_{\red}|>2\eps N$, by~\eqref{eq:complete-codegree} there is a vertex~$a\in V'_{\red}$ such that~$z_2z_3a$ and~$az_3z_4$ are blue edges in~$E^*$. Hence, $z_1z_2z_3$ and~$z_3z_4z_5$ are in the same blue component via~$z_2z_3a$ and~$az_3z_4$. Finally, consider two disjoint blue edges~$z_1z_2z_3$ and~$z_4z_5z_6$ in~$E^*$. Since~$|V'_{\red}|>2\eps N$, by~\eqref{eq:complete-codegree} there is a vertex~$a\in V'_{\red}$ such that~$z_2z_3a$ and~$az_4z_5\in E^*$ are blue edges in~$E^*$. Since~$z_2z_3a$ and~$az_4z_5$ are blue edges in~$E^*$ which share one vertex, we conclude that~$z_1z_2z_3$ and~$z_4z_5z_6$ are in the same blue component via~$z_2z_3a$ and~$az_4z_5$.
    \end{claimproof}

    Clearly, the colours red and blue may be interchanged in Claims~\ref{claim:all-blue-monoxcomp} and~\ref{claim:pair-few-blue-cross-edges}; in particular, there is a red component~$C'_{\red}$ such that~$E(H)\cap\binom{V'_{\blue}}{3}\subseteq C'_{\red}$. If~$V'_{\blue}$ contains at least~$3s+3\eps N$ vertices, then by~\eqref{eq:complete-codegree} and Lemma~\ref{lem:setup-properties}\ref{item:complete-active} we may greedily find both a copy of~$M^{(3)}_s$ and a copy of~$C^{(3)}_4$ in~$C'_{\red}$. If not, by Lemma~\ref{lem:setup-properties}\ref{item:colour-bipartition-near-complete}~$V'_{\red}$ contains at least~$3t+3\eps N$ vertices, so by~\eqref{eq:complete-codegree} and Lemma~\ref{lem:setup-properties}\ref{item:complete-active} we may greedily find both a copy of~$M^{(3)}_t$ and a copy of~$C^{(3)}_4$ in~$C'_{\blue}$; this gives a contradiction.

    \textbf{Case 2:}~$|V'_{\red}|<s'$ or~$|V'_{\blue}|<s'$. We shall assume~$|V'_{\red}|<s'$; the proof for~$|V'_{\blue}|<s'$ will be identical with red and blue interchanged. We first show that~$A$ contains a connected matching of size~$s'$ whose vertex set is contained in~$V'_{\blue}$.

    \begin{claim} \label{claim:majority-colour-matching-dominant-monoxcomp}
        There is a matching~$M_A$ of size~$s'$ in~$A$ with~$V(M_A)\subseteq V'_{\blue}$.
    \end{claim}

    \begin{claimproof}
        By Lemma~\ref{lem:setup-properties}\ref{item:colour-bipartition-near-complete} we have~$|V'_{\blue}| \ge (4+\eta-2\delta)s \ge (4+9\delta)s$, so by~\eqref{eq:xdiv-xdiv-contra} and Lemma~\ref{lem:matching-div-div} applied to the induced colouring~$H[V'_{\blue}]=H_{\red}[V'_{\blue}]\cup H_{\blue}[V'_{\blue}]$ we obtain a matching~$M'$ of size~$s(1+2\delta)$ in a monochromatic component~$C$ of~$H$ with~$C^{(3)}_4\nsubseteq C$ and~$V(M')\subseteq V'_{\blue}$; note that~$|V'_{\blue}\sm V(M')|\ge s(1+2\delta)$. If~$C$ is blue, then by~\eqref{eq:xdiv-xdiv-contra} and Lemma~\ref{lem:not-scarlet-matching-force-blue-C4-matching} applied to~$M'$ with the colours red and blue interchanged, we obtain~$C=A$; in this case, we simply take~$M_A=M'$. If~$C$ is red, Lemma~\ref{lem:matching-pivot-no-C4} applied to~$M'$ and~$V'_{\blue}\sm V(M')$ returns the desired connected matching~$M_A$ of size at least~$s(1+2\delta)-5\eps N \ge s'$ in~$A$ with~$V(M_A)\subseteq V'_{\blue}$.
    \end{claimproof}

    By Claim~\ref{claim:majority-colour-matching-dominant-monoxcomp} and~\eqref{eq:xdiv-xdiv-contra} we have~$C^{(3)}_4\nsubseteq A$. Since~$\frac{N}{5}(1-26\eps^{1/12})\ge s'>|V'_{\red}|$, Lemma~\ref{lem:build-blue-C4-matching-scarlet-C5} applied to~$M_A$ with the colours red and blue interchanged and~\eqref{eq:xdiv-xdiv-contra} implies that~$C^{(3)}_5\subseteq A$. But now~$A$ contains both~$M_A$ and~$C^{(3)}_5$, which contradicts~\eqref{eq:xdiv-xdiv-contra}.
\end{proof}

\begin{proof}[Proof of Lemma~\ref{lem:matching-div-xdiv}]
    By Lemmas~\ref{lem:vertex-shadow-monoxcomp} and~\ref{lem:dominant-monoxcomp}, we may assume Setup~\ref{setup} with~$N=(5+\eta)s$ from now onwards. Suppose that the~$2$-edge-colouring~$H=H_{\red}\cup H_{\blue}$ is such that
    \begin{equation} \label{eq:div-xdiv-contra}
        \begin{split}
            &\textrm{no red component contains }\;M^{(3)}_s\;\textrm{, and } \\
            &\textrm{no blue component contains }\;M^{(3)}_t\;\textrm{ and either }\;C^{(3)}_4\;\textrm{ or }\;C^{(3)}_5\,.
        \end{split}
    \end{equation}
    We distinguish two cases and obtain a contradiction to~\eqref{eq:div-xdiv-contra} in each of them. Set~$\delta:=\frac{\eta}{20}$.

    \textbf{Case 1:}~$|V'_{\red}|\ge s':=(1+\delta)s$. Partition~$V'=U\cup U'$ with~$U'\subseteq V'_{\red}$, $|U'|=s'$ and~$U=V'\sm U'$. By Lemma~\ref{lem:setup-properties}\ref{item:colour-bipartition-near-complete} we have~$|U| \ge (4+\eta-2\delta)s \ge (4+4\delta)s$, so by~\eqref{eq:div-xdiv-contra} and Lemma~\ref{lem:matching-div-div} applied to the induced colouring~$H[U]=H_{\red}[U]\cup H_{\blue}[U]$ we obtain a matching~$M'$ of size~$s'$ in a blue component~$C$ of~$H$ with~$C^{(3)}_4\nsubseteq C$ and~$V(M')\subseteq U$. Now Lemma~\ref{lem:matching-pivot-no-C4} applied to~$M'$ and~$U'$ with the colours red and blue interchanged returns a connected matching~$M$ of size at least~$s'-5\eps N \ge s$ in~$S$, which contradicts~\eqref{eq:div-xdiv-contra}.

    \textbf{Case 2:}~$|V'_{\red}|<s'$. By Lemma~\ref{lem:setup-properties}\ref{item:colour-bipartition-near-complete} we have~$|V'_{\blue}| \ge (4+\eta-2\delta)s \ge (4+4\delta)s$, so by~\eqref{eq:div-xdiv-contra} and Lemma~\ref{lem:matching-div-div} applied to the induced colouring~$H[V'_{\blue}]=H_{\red}[V'_{\blue}]\cup H_{\blue}[V'_{\blue}]$ we obtain a matching~$M'$ of size~$s'$ in a blue component~$C$ of~$H$ with~$C^{(3)}_4\nsubseteq C$ and~$V(M')\subseteq V'_{\blue}$. Now by~\eqref{eq:div-xdiv-contra} and Lemma~\ref{lem:not-scarlet-matching-force-blue-C4-matching} applied to~$M'$ with the colours red and blue interchanged, we obtain~$C=A$. Since~$\frac{N}{5}(1-26\eps^{1/12})\ge s'>|V'_{\red}|$, Lemma~\ref{lem:build-blue-C4-matching-scarlet-C5} applied to~$M'$ with the colours red and blue interchanged and~\eqref{eq:div-xdiv-contra} yields~$C^{(3)}_5\subseteq A$. But now~$A$ contains both~$M'$ and~$C^{(3)}_5$, which contradicts~\eqref{eq:div-xdiv-contra}.
\end{proof}

\begin{proof}[Proof of Lemma~\ref{lem:matching-xdiv-div}]
    By Lemmas~\ref{lem:vertex-shadow-monoxcomp} and~\ref{lem:dominant-monoxcomp}, from now onwards we may assume Setup~\ref{setup} with~$N=\max(3s+t,5t)+\eta s$. Suppose that the~$2$-edge-colouring~$H=H_{\red}\cup H_{\blue}$ is such that
    \begin{equation} \label{eq:xdiv-div-contra}
        \begin{split}
            &\textrm{no red component contains }\;M^{(3)}_s\;\textrm{ and either }\;C^{(3)}_4\;\textrm{ or }\;C^{(3)}_5\;\textrm{, and } \\
            &\textrm{no blue component contains }\;M^{(3)}_t\,.
        \end{split}
    \end{equation}
    We distinguish three cases and obtain a contradiction to~\eqref{eq:xdiv-div-contra} in each of them. Set~$\delta:=\frac{\eta}{20}$, $t':=t+\delta s$ and~$s'=s(1+\delta)$.

    \textbf{Case 1:}~$|V'_{\blue}|\ge 2t'$. By Lemma~\ref{lem:matching-div-div} applied to the colouring~$H=H_{\red}\cup H_{\blue}$ and~\eqref{eq:xdiv-div-contra}, there is a matching of size~$\frac{N-t}{3+\delta} \ge t'$ in a red component~$C$ of~$H$ with~$C^{(3)}_4\nsubseteq C$; let~$M$ be a matching of size~$t'$ in~$C$ with the maximum number of vertices in~$V'_{\red}$ among all matchings of size~$t'$ in~$C$. Set~$W:=V'_{\red}\sm V(M)$, $X:=V'_{\blue}\sm V(M)$ and~$M^*:=\{e\in M:|e\cap V'_{\blue}|\ge2\}$; note that
    \begin{equation} \label{eq:M^*-size-red}
        |M^*|\ge\min\left(|M|,\max\left(\frac{|V'_{\blue}|-|M|-|X|}{2},0\right)\right)\,.
    \end{equation}
    Lemma~\ref{lem:matching-pivot-max-red-vertices} applied to~$M$ and~$W$ returns a matching~$M'$ of size~$\min(|M^*|,\max(|W|-2\eps N,0))$ in~$A$ with a function~$\phi'\colon M'\to M^*$. Set~$M_1:=M\sm\phi'(M')$ and~$X_1:=V'_{\blue}\sm(V(M')\cup V(M_1))$; note that~$|X_1|=|X|+|M'|$. Now Lemma~\ref{lem:matching-pivot-no-C4} applied to~$M_1$ and~$X_1$ returns a matching~$M''$ of size~$\min(|M_1|,|X_1|) - 5\eps N$ in~$A$ with a function~$\phi''\colon M''\to M_1$. The properties of~$\phi'$ and~$\phi''$ imply that~$M'\cup M''$ is a connected matching in~$A$ of size
    \begin{align*}
        |M'| + \min(|M_1|,|X_1|) - 5\eps N & = \min(|M|,2|M'|+|X|) - 5\eps N \\
            & \ge \min(|M|,|V'_{\blue}|-|M|,N(1-2\eps)-|V(M)|) - 5\eps N \ge t\,,
    \end{align*}
    where the inequalities are obtained by applying~\eqref{eq:M^*-size-red}, the assumption that~$|V'_{\blue}|\ge 2t'$ and Lemma~\ref{lem:setup-properties}\ref{item:colour-bipartition-near-complete}. But this contradicts~\eqref{eq:xdiv-div-contra}.

    \textbf{Case 2:}~$t'\le|V'_{\blue}|<2t'$. Fix~$U''\subseteq U'\subseteq V'_{\blue}$ with~$|U'|=t'$ and~$|U''|=\max(4t-3s',0)$, and set~$U:=V'\sm U'$ and~$U^*:=V'\sm U''$. By Lemma~\ref{lem:setup-properties}\ref{item:colour-bipartition-near-complete} we have~$|U| \ge 4t+(\eta-2\delta)s \ge 4t+4\delta s$, so by~\eqref{eq:xdiv-div-contra} and Lemma~\ref{lem:matching-div-div} applied to the induced colouring~$H[U]=H_{\red}[U]\cup H_{\blue}[U]$ we obtain a matching~$M$ of size~$t'$ in a red component~$C$ of~$H$ with~$V(M)\subseteq U$; Lemma~\ref{lem:not-scarlet-matching-force-blue-C4-matching} applied to~$M$ with~\eqref{eq:xdiv-div-contra} yields~$C=S$. By~\eqref{eq:xdiv-div-contra} and Lemma~\ref{lem:matching-div-div} applied to the induced colouring~$H[U^*]=H_{\red}[U^*]\cup H_{\blue}[U^*]$, there is a matching~$M^*$ of size~$s'$ in a red component~$C^*$ of~$H$ with~$C^{(3)}_4\nsubseteq C^*$. For~$i\in[3]\cup\{0\}$ set~$P^*_i:=\{e\in M^*:|e\cap V'_{\red}|=i\}$. Since~$|V'_{\blue}|<2t'$ and~$U''\subseteq V'_{\blue}$, we obtain~$|V'_{\blue}|-|U''| \ge |V'_{\blue}\cap V(M^*)| \ge 3|P^*_0| + 2|P^*_1| + |P^*_2| \ge \frac{3s'}{2} + |P^*_0|-|P^*_2|-|P^*_3|$, which yields
    \begin{equation} \label{eq:bound-red-imbalance-case1}
        |P^*_0|-|P^*_2|-|P^*_3| < 2t' - 3s'/2 - \max(4t-3s',0)\,.
    \end{equation}
    Suppose that~$C^*\ne S$. Then Lemma~\ref{lem:not-scarlet-matching-force-blue-C4-matching} applied to~$M^*$ with~\eqref{eq:bound-red-imbalance-case1} returns a blue component containing a copy of~$C^{(3)}_4$ and a matching of size
    \[s'\left(5/2-24\eps^{1/2}\right) - 2t' + \max(4t-3s',0) - 2\eps^{1/2}N - 2 \ge t\,,\]
    which contradicts~\eqref{eq:xdiv-div-contra}. Hence, we have~$C^*=S$ and therefore~$C^{(3)}_4\nsubseteq S$. Now Lemma~\ref{lem:matching-pivot-no-C4} applied to~$M'$ and~$U'$ returns a connected matching of size at least~$t'-5\eps N \ge t$ in~$A$, which contradicts~\eqref{eq:xdiv-div-contra}.

    \textbf{Case 3:}~$|V'_{\blue}|<t'$. By Lemma~\ref{lem:setup-properties}\ref{item:colour-bipartition-near-complete} we have~$|V'_{\red}| \ge 4t+(\eta-2\delta)s \ge 4t+4\delta s$, so by~\eqref{eq:xdiv-div-contra} and Lemma~\ref{lem:matching-div-div} applied to the induced colouring~$H[V'_{\red}]=H_{\red}[V'_{\red}]\cup H_{\blue}[V'_{\red}]$ we obtain a matching~$M$ of size~$t'$ in a red component~$C$ of~$H$ with~$V(M)\subseteq V'_{\red}$; Lemma~\ref{lem:not-scarlet-matching-force-blue-C4-matching} applied to~$M$ with~\eqref{eq:xdiv-div-contra} yields~$C=S$. By~\eqref{eq:xdiv-div-contra} and Lemma~\ref{lem:matching-div-div} applied to the colouring~$H=H_{\red}\cup H_{\blue}$, there is a matching~$M^*$ of size~$s'$ in a red component~$C^*$ of~$H$ with~$C^{(3)}_4\nsubseteq C^*$. For~$i\in[3]\cup\{0\}$ set~$P^*_i:=\{e\in M^*:|e\cap V'_{\red}|=i\}$. Since~$|V'_{\blue}|<t'$, we obtain~$|V'_{\blue}| \ge |V'_{\blue}\cap V(M^*)| \ge 3|P^*_0| + 2|P^*_1| + |P^*_2| \ge \frac{3s'}{2} + |P^*_0|-|P^*_2|-|P^*_3|$, which yields
    \begin{equation} \label{eq:bound-red-imbalance-case2}
        |P^*_0|-|P^*_2|-|P^*_3| < t' - 3s'/2\,.
    \end{equation}
    Suppose that~$C^*\ne S$. Then Lemma~\ref{lem:not-scarlet-matching-force-blue-C4-matching} applied to~$M^*$ with~\eqref{eq:bound-red-imbalance-case2} returns a blue component containing a copy of~$C^{(3)}_4$ and a matching of size
    \[s'\left(5/2-24\eps^{1/2}\right) - t' - 2\eps^{1/2}N - 2 \ge t\,,\]
    which contradicts~\eqref{eq:xdiv-div-contra}. Hence, we have~$C^*=S$ and therefore~$C^{(3)}_4\nsubseteq S$. Now Lemma~\ref{lem:build-blue-C4-matching-scarlet-C5} applied to~$M$ combined with~$\frac{N}{5}(1-26\eps^{1/12})\ge t$ and~\eqref{eq:xdiv-div-contra} yields~$C^{(3)}_5\subseteq S$. But this means that~$S$ contains both~$M^*$ and~$C^{(3)}_5$, contradicting~\eqref{eq:xdiv-div-contra}.
\end{proof}

\section{Hypergraph regularity and proving Lemma~\ref{lem:monox-cycle-embedding}} \label{sec:regularity-cycle-embedding}

In this section, we prove Lemma~\ref{lem:monox-cycle-embedding}. We first introduce the formulation of hypergraph regularity we require, closely following the presentation of Allen, B\"{o}ttcher, Cooley and Mycroft~\cite{AllenBoettcherCooleyMycroft}.

Let~$\cH$ be a hypergraph. We say that~$\cH$ is a \emph{complex} if its edge set~$E(\cH)$ is \emph{down-closed}: if~$e\in E(\cH)$ and~$f\subseteq e$, then~$f\in E(\cH)$; it is a~\emph{$k$-complex} if it is a complex whose edges are all of size at most~$k$. Write~$\cH^{(i)}$ for the~$i$-graph on~$V(\cH)$ formed of the edges of size~$i$ of~$\cH$. For a subset~$U\subseteq V(\cH)$, the \emph{induced subgraph}~$\cH[U]$ is the hypergraph with vertex set~$U$ and edge set~$\{e\in E(\cH):e\subseteq U\}$. Let~$\cP$ be a partition of~$V(\cH)$ into parts~$V_1,\dots,V_s$. We say that a subset~$S\subseteq V(\cH)$ is \emph{$\cP$-partite} for all~$i\in[s]$ we have~$|S\cap V_i|\le1$. We say that~$\cH$ is~$\cP$-partite if all its edges are~$\cP$-partite; in this case, we call the parts of~$\cP$ the \emph{vertex classes} of~$\cH$. We say that~$\cH$ is~$s$-partite if there is a partition~$\cP$ of~$V(\cH)$ into~$s$ parts for which~$\cH$ is~$\cP$-partite.

Let~$\cH$ be a~$\cP$-partite hypergraph and~$X$ be a~$j$-set of clusters of~$\cH$. Write~$\cH_X$ for the~$j$-partite~$j$-uniform subgraph of~$\cH^{(j)}$ induced on~$\bigcup X$ whose vertex classes are the elements of~$X$. Write~$\cH_{X^<}$ for the~$j$-partite hypergraph with vertex set~$\bigcup X$ and edge set~$\bigcup_{X'\subsetneq X}\cH_{X'}$.

Let~$i\ge2$ be an integer, and let~$\cP$ be a partition of a vertex set~$V$ into~$i$ parts. Let~$H_i$ and~$H_{i-1}$ be a~$\cP$-partite~$i$-graph and a~$\cP$-partite~$(i-1)$-graph, respectively, on a common vertex set~$V$. Let~$K_i(H_{i-1})$ denote the~$\cP$-partite~$i$-graph on~$V$ whose edges are all~$\cP$-partite~$i$-sets in~$V$ which are \emph{supported} on~$H_{i-1}$, that is, all~$\cP$-partite~$i$-sets in~$V$ which induce a copy of the complete~$(i-1)$-graph~$K^{(i-1)}_i$ on~$i$ vertices in~$H_{i-1}$. The \emph{density of~$H_i$ with respect to~$H_{i-1}$} is defined to be
\[d(H_i\mid H_{i-1}):=\frac{|K_i(H_{i-1})\cap H_i|}{|K_i(H_{i-1})|}\]
if~$|K_i(H_{i-1})|>0$ and~$d(H_i\mid H_{i-1}):=0$ if~$|K_i(H_{i-1})|=0$. More generally, for a collection~$\bQ=(Q_1,\dots,Q_r)$ of~$r$ (not necessarily disjoint) subgraphs of~$H_{i-1}$, we define~$K_i(\bQ):=\bigcup_{j\in[r]}K_i(Q_j)$, 
\[d(H_i\mid\bQ):=\frac{|K_i(\bQ)\cap H_i|}{|K_i(\bQ)|}\]
if~$|K_i(\bQ)|>0$ and~$d(H_i\mid \bQ):=0$ if~$|K_i(\bQ)|=0$. We say that~$H_i$ is~\emph{$(d,\eps,r)$-regular} with respect to~$H_{i-1}$ if for every~$r$-set~$\bQ$ of subgraphs of~$H_{i-1}$ with~$|K_i(\bQ)|>\eps|K_i(H_{i-1})|$ we have~$d(H_i\mid\bQ)=d\pm\eps$. We say that~$H_i$ is~\emph{$(\eps,r)$-regular} with respect to~$H_{i-1}$ if there exists some~$d$ for which~$H_i$ is~$(d,\eps,r)$-regular with respect to~$H_{i-1}$. Given an~$i$-graph~$G$ whose vertex set contains that of~$H_{i-1}$, we say that~$G$ is~\emph{$(d,\eps,r)$-regular} with respect to~$H_{i-1}$ if the~$i$-partite subgraph of~$G$ induced by the vertex classes of~$H_{i-1}$ is~\emph{$(d,\eps,r)$-regular} with respect to~$H_{i-1}$. The \emph{density of~$G$ with respect to~$H_{i-1}$} is the density of the~$i$-partite subgraph of~$G$ induced by the vertex classes of~$H_{i-1}$ with respect to~$H_{i-1}$.

Let~$s\ge k\ge3$ be integers, and let~$\cH$ be an~$s$-partite~$k$-complex on vertex classes~$V_1,\dots,V_s$. For any~$A\subseteq[s]$ write~$V_A$ for~$\bigcup_{i\in A}V_i$. We say that~$\cH$ is~\emph{$(d_k,\dots,d_2,\eps_k,\eps,r)$-regular} if
\begin{enumerate}[label=(\alph{*})]
    \item for any~$2\le i\le k-1$ and any~$A\in\binom{[s]}{i}$ the induced subgraph~$\cH^{(i)}[V_A]$ is~$(d_i,\eps,1)$-regular with respect to~$\cH^{(i-1)}[V_A]$, and
    \item for any~$A\in\binom{[s]}{k}$ the induced subgraph~$\cH^{(k)}[V_A]$ is~$(d_k,\eps_k,r)$-regular with respect to~$\cH^{(k-1)}[V_A]$.
\end{enumerate}
For a~$(k-1)$-tuple~$\bd=(d_k,\dots,d_2)$ we write~$(\bd,\eps_k,\eps,r)$-regular to mean~$(d_k,\dots,d_2,\eps_k,\eps,r)$-regular. We say that a~$(k-1)$-complex~$\cJ$ is~\emph{$(t_0,t_1,\eps)$-equitable} if it has the following properties.
\begin{enumerate}[label=(\alph{*})]
    \item There is a partition~$\cP$ of~$V(\cJ)$ into~$t$ parts, where~$t_0\le t\le t_1$, of equal size such that~$\cJ$ is~$\cP$-partite; we call~$\cP$ the \emph{ground partition} of~$\cJ$ and the parts of~$\cP$ the \emph{clusters} of~$\cJ$.
    \item There is a \emph{density vector}~$\bd=(d_{k-1},\dots,d_2)$ such that~$\cJ$ is~$(\bd,\eps,\eps,1)$-regular and for each~$2\le i\le k-1$ we have~$d_i\ge1/t_1$ and~$1/d_i\in\NN$.
\end{enumerate}
For any~$k$-set~$X$ of clusters of~$\cJ$, we write~$\Jchat_X$ for the~$k$-partite~$(k-1)$-graph~$\cJ_{X^ {<}}^{(k-1)}$ and call it a \emph{polyad}. Given a~$(t_0,t_1,\eps)$-equitable~$(k-1)$-complex~$\cJ$, a~$k$-graph~$G$ on~$V(\cJ)$ and a~$k$-set~$X$ of clusters of~$\cJ$, we say that~$G$ is~\emph{$(\eps_k,r)$-regular with respect to~$X$} if there exists some~$d$ such that~$G$ is~$(d,\eps_k,r)$-regular with respect to the polyad~$\Jchat_X$. We write~$d^*_{G,\cJ}(X)$ for the density of~$G$ with respect to~$\Jchat_X$; we often omit~$G$ and~$\cJ$ from the subscripts when they are clear from the context.

We give the key definition of a regular slice.

\begin{definition}[Regular slice] \label{def:regular-slice}
    Given positive integers~$k,r,t_0,t_1\in\NN$, positive constants~$\eps,\eps_k>0$, a~$k$-graph~$G$ and a~$(k-1)$-complex~$\cJ$ on~$V(G)$, we call~$\cJ$ a~\emph{$(t_0,t_1,\eps,\eps_k,r)$-regular slice} for~$G$ if~$\cJ$ is~$(t_0,t_1,\eps)$-equitable and~$G$ is~$(\eps_k,r)$-regular with respect to all but at most~$\eps_k\binom{t}{k}$ of the~$k$-sets of clusters of~$\cJ$, where~$t$ is the number of clusters of~$\cJ$.
\end{definition}

It will sometimes be convenient to not specify the parameters~$t_0$ and~$t_1$ when we instead specify the number of clusters~$t$ and the density vector~$\bd$; in these situations, we write that~$\cJ$ is~$(\cdot,\cdot,\eps)$-equitable or that~$\cJ$ is a~$(\cdot,\cdot,\eps,\eps_k,r)$-regular slice for~$G$.

We state the version of the Regular Slice Lemma which we need, which is the special case of~\cite[Lemma 10]{AllenBoettcherCooleyMycroft} given in~\cite[Lemma 11]{LoPfenninger2023}.

\begin{lemma}[Regular Slice Lemma] \label{lem:regular-slice}
    Given an integer~$k\ge3$, for all~$t_0\in\NN$, $\eps_k>0$ and functions~$r\colon\NN\to\NN$ and~$\eps\colon\NN\to[0,1]$, there are positive integers~$t_1$ and~$n_0$ such that the following holds for all~$n\ge n_0$ which are divisible by~$t_1!$. Suppose that~$K=K_{\red}\cup K_{\blue}$ is a~$2$-edge-coloured complete~$k$-graph on~$n$ vertices. Then there exists a~$(k-1)$-complex~$\cJ$ on~$V(K)$ which is a~$(t_0,t_1,\eps(t_1),\eps_k,r(t_1))$-regular slice for both~$K_{\red}$ and~$K_{\blue}$.
\end{lemma}

\begin{definition}[The~$d$-reduced~$k$-graph]
    Let~$G$ be a~$k$-graph and~$\cJ$ be a~$(t_0,t_1,\eps,\eps_k,r)$-regular slice for~$G$. For~$d>0$ the~\emph{$d$-reduced~$k$-graph}~$R_d(G)$ is the~$k$-graph whose vertices are the clusters of~$\cJ$ and whose edges are all~$k$-sets~$X$ of clusters of~$\cJ$ such that~$G$ is~$(\eps_k,r)$-regular with respect to~$X$ and~$d^*(X)\ge d$.
\end{definition}

A \emph{fractional matching} in a~$k$-graph~$H$ is a function~$w\colon E(H)\to[0,1]$ such that for every~$v\in V(H)$ we have~$\sum_{e\in E(H):v\in e}w(e)\le1$; its \emph{weight} is~$\sum_{e\in E(H)}w(e)$, and it is \emph{tightly connected} if all the edges in~$\{e\in E(H):w(e)>0\}$ lie in a single tight component.

\begin{lemma}[{\cite[Lemma~13]{AllenBoettcherCooleyMycroft}}] \label{lem:cycle-embedding-divisible}
    Given~$k,r,n_0,t\in\NN$ and~$\psi,\eps,\eps_k,d_2,\dots,d_k>0$ such that~$1/d_i\in\NN$ for each~$2\le i\le k-1$ and
    \[\frac{1}{n_0}\ll\frac{1}{t}\;, \frac{1}{n_0}\ll\frac{1}{r},\eps\ll\eps_k,d_2,\dots,d_{k-1}\;\textrm{ and }\eps_k\ll\psi,d_k,\frac{1}{k}\,, \]
    the following holds for all integers~$n\ge n_0$. Let~$G$ be a~$k$-graph on~$n$ vertices and let~$\cJ$ be a~$(\cdot,\cdot,\eps,\eps_k,r)$-regular slice for~$G$ with~$t$ clusters and density vector~$\bd=(d_{k-1},\dots,d_2)$. Suppose that~$R_{d_k}(G)$ contains a tightly connected fractional matching with weight~$\mu$. Then~$G$ contains a tight cycle of length~$\ell$ for any~$\ell\le(1-\psi)\mu kn/t$ divisible by~$k$.
\end{lemma}

We also need the following version of Lemma~\ref{lem:cycle-embedding-divisible} to obtain cycles of lengths not divisible by~$k$. The proof is mostly identical to that of~\cite[Lemma~13]{AllenBoettcherCooleyMycroft}, so we omit a detailed proof. Indeed, the main modification needed is an extension of the tight walk in~$R_{d_k}$ which serves as the `frame' for the construction of the desired tight cycle so that it `visits' an edge of the given copy of~$C^{(k)}_p$; the congruence class modulo~$k$ is adjusted by `going around' the copy of~$C^{(k)}_p$ suitably many times. For sufficiently short cycles, it is enough to walk in the given copy of~$C^{(k)}_p$. The lower bound on~$\ell$ derives from the fact that non-divisible lengths require up to~$k-1$ remainder adjustments each of length at most~$p(k+1)\le(2k-1)(k+1)$; we have~$(2k-1)(k+1)(k-1)=(2k^2-k-2)k+1$.

\begin{lemma} \label{lem:cycle-embedding-all}
    Given positive integers~$k,r,n_0,t\in\NN$ and positive constants~$\psi,\eps,\eps_k,d_2,\dots,d_k$ such that~$1/d_i\in\NN$ for each~$2\le i\le k-1$ and
    \[\frac{1}{n_0}\ll\frac{1}{t}\;, \frac{1}{n_0}\ll\frac{1}{r},\eps\ll\eps_k,d_2,\dots,d_{k-1}\;\textrm{ and }\eps_k\ll\psi,d_k,\frac{1}{k}\,, \]
    the following holds for all integers~$n\ge n_0$. Let~$G$ be a~$k$-graph on~$n$ vertices and let~$\cJ$ be a~$(\cdot,\cdot,\eps,\eps_k,r)$-regular slice for~$G$ with~$t$ clusters and density vector~$\bd=(d_{k-1},\dots,d_2)$. Suppose that~$R_{d_k}(G)$ has a tight component which contains both a tightly connected fractional matching with weight~$\mu$ and a copy of~$C^{(k)}_p$ for some~$p\in\NN$ satisfying~$k<p<2k$ and~$\gcd(k,p)=1$. Then~$G$ contains a tight cycle of length~$\ell$ for any~$(2k^2-k-2)k\le\ell\le(1-\psi)\mu kn/t$.
\end{lemma}

The following lemma strengthens the almost-completeness properties of a~$3$-graph by deleting a small fraction of its vertices. A (fairly standard) proof can be found in~\cite[Lemma 4.1]{HaxellLuczakPengRodl2006}.

\begin{lemma}[{\cite[Lemma 8.1]{HaxellLuczakPengRodl}}] \label{lem:subhypergraph-active-large-codegree}
    For all~$\eps>0$ and~$n\in\NN$ every~$3$-graph on~$n$ vertices with at least~$(1-\eps)\binom{n}{3}$ edges contains a~$10\eps^{1/6}$-complete subhypergraph with~$n'\ge(1-10\eps^{1/6})n$ vertices.
\end{lemma}

We are now ready to prove Lemma~\ref{lem:monox-cycle-embedding}.

\begin{proof}[Proof of Lemma~\ref{lem:monox-cycle-embedding}]
    Let~$d_3=\frac{1}{2}$ and~$\frac{1}{t_0}\ll\eps_3\ll\delta$. Pick functions~$r(\cdot)$ growing sufficiently quickly and~$\eps(\cdot)$ tending to zero sufficiently quickly such that for any~$t'\ge t_0$ and~$d_2\ge\frac{1}{t'}$ we may apply Lemmas~\ref{lem:cycle-embedding-divisible} and~\ref{lem:cycle-embedding-all} with~$r(t')$ and~$\eps(t')$. Lemma~\ref{lem:regular-slice} with inputs~$t_0$, $\eps_3$, $r(\cdot)$ and~$\eps(\cdot)$ returns positive integers~$t_1$ and~$n_0$. Set~$\eps=\eps(t_1)$ and~$r=r(t_1)$. Pick~$n_1$ large enough to apply Lemmas~\ref{lem:cycle-embedding-divisible} and~\ref{lem:cycle-embedding-all} for all~$n\ge n_1$ and set~$n_2:=n_1+t_1!$.

    Suppose that~$K=K_{\red}\cup K_{\blue}$ be a~$2$-edge-coloured complete~$3$-graph on~$n\ge n_2$ vertices. Let~$K'=K'_{\red}\cup K'_{\blue}$ be a complete subgraph of~$K$ on~$n':=t_1!\cdot\left\lfloor\frac{n}{t_1!}\right\rfloor\ge n_1$ vertices with the induced~$2$-edge-colouring. Then Lemma~\ref{lem:regular-slice} returns a~$(t_0,t_1,\eps,\eps_3,r)$-regular slice~$\cJ$ for both~$K'_{\red}$ and~$K'_{\blue}$. Let~$t$ be the number of clusters of~$\cJ$ and~$(d_2)$ be the density vector of~$\cJ$. Let~$H=R_{d_3}(K'_{\red})\cup R_{d_3}(K'_{\blue})$ be a~$3$-graph with a~$2$-edge-colouring~$H=H_{\red}\cup H_{\blue}$ such that~$R_{d_3}(K'_{\red})\sm R_{d_3}(K'_{\blue})\subseteq H_{\red}$ and~$R_{d_3}(K'_{\blue})\sm R_{d_3}(K'_{\red})\subseteq H_{\blue}$.

    \begin{claim} \label{claim:reduced-graph-almost-complete}
        $H$ has at least~$(1-2\eps_3)\binom{t}{3}$ edges.
    \end{claim}

    \begin{claimproof}
        Since~$\cJ$ is a~$(\cdot,\cdot,\eps,\eps_3,r)$-regular slice for both~$K'_{\red}$ and~$K'_{\blue}$, at least~$(1-2\eps_3)\binom{t}{3}$ triples~$X$ of clusters of~$\cJ$ are such that both~$K'_{\red}$ and~$K'_{\blue}$ are~$(\eps_3,r)$-regular with respect to~$X$. Let~$X$ be such a~triple. Since~$K'_{\red}$ and~$K'_{\blue}$ are complements of each other, we have~$d^*_{K'_{\red}}(X) + d^*_{K'_{\blue}}(X) = 1$. Hence, we have~$d^*_{K'_{\red}}(X)\ge d_3=\frac{1}{2}$ or~$d^*_{K'_{\blue}}(X)\ge d_3=\frac{1}{2}$, which implies that~$X$ is an edge of~$R_{d_3}(K'_{\red})$ or~$R_{d_3}(K'_{\blue})$.
    \end{claimproof}

    Set~$\eps':=20\eps_3^{1/6}$. Claim~\ref{claim:reduced-graph-almost-complete} and Lemma~\ref{lem:subhypergraph-active-large-codegree} applied to~$H$ gives an~$\eps'$-complete subhypergraph~$H'$ on~$t'\ge(1-\eps')t$ vertices. Since~$\eps'\ll\delta^2$, for all~$i\in[4]$ and~$\beta\in(0,1]$ we have~$m^*_i(\beta,\eps',t')\ge f(\beta)(m^*_i(f)-\delta^2)t$. To complete the proof for~\ref{item:monox-cycle-embedding-div-div}, we observe that the definition of~$m^*_i(\beta,\eps',t')$ implies that~$H'$ contains a matching of size at least~$m^*_1(\beta,\eps',t')$ whose edges all lie in a single red tight component or a matching of size at least~$\beta m^*_1(\beta,\eps',t')$ whose edges all lie in a single blue tight component. By applying Lemma~\ref{lem:cycle-embedding-divisible} to the monochromatic matching mentioned above and noting that~$(1-\delta^2)f(\beta)(m^*_i(f)-\delta^2)3n'\ge(f(\beta)m^*_i(f)-\delta)3n$, we obtain the desired outcome. Finally, we remark that the proofs of~\ref{item:monox-cycle-embedding-xdiv-xdiv}--\ref{item:monox-cycle-embedding-xdiv-div} are entirely analogous to that of~\ref{item:monox-cycle-embedding-div-div}; indeed, the main difference is that we apply Lemma~\ref{lem:cycle-embedding-all} instead of Lemma~\ref{lem:cycle-embedding-divisible} whenever the monochromatic tight component containing the edges of the relevant matching also contains a copy of~$C^{(3)}_4$ or~$C^{(3)}_5$.
\end{proof}

\section{Concluding remarks}

In this paper, we asymptotically determined the off-diagonal Ramsey numbers of~$3$-uniform tight cycles. A keen reader may notice that the statement of Theorem~\ref{thm:main}\ref{item:thm-upper-bound} contains a somewhat unusual lower bound of~$t\ge13$. The underlying technical reason is that the unavoidable walk-traversing cycle-extending steps for non-divisible cycle lengths in the proof of Lemma~\ref{lem:cycle-embedding-all} (which is mostly identical to that of~\cite[Lemma~13]{AllenBoettcherCooleyMycroft}) are more `expensive' than ($k$-divisible) in-situ cycle-extending steps. We anticipate that it may be possible to obtain shorter cycles of non-divisible lengths with further technical analysis or different residue-switching gadgets.

A natural follow-up direction is the corresponding off-diagonal problem for general~$k\ge4$. Conjecture~\ref{conj:k-cycle-ramsey} for the diagonal case suggests that~$\gcd(k,i)$ is the key parameter which generalises the role of divisibility by~$3$, and the structure of Theorem~\ref{thm:main} indicates that the phase diagram should become substantially richer. Progress in this direction will require the development of substantial new theory to obtain cycles of non-divisible lengths (even for the diagonal case, which remains unresolved) and the development of new sharp extremal constructions.

\section*{Acknowledgements}

The author thanks Hong Liu for stimulating discussions and insightful comments.


\bibliography{references} 
\bibliographystyle{abbrv}


\appendix

\section{Proving Lemma~\ref{lem:matching-div-div}} \label{sec:matching-div-div}

In this appendix, we prove Lemma~\ref{lem:matching-div-div}, which is a generalisation of~\cite[Lemma~2.1]{HaxellLuczakPengRodl}. The proof follows by observing that the techniques and arguments used in~\cite{HaxellLuczakPengRodl} immediately imply the more general conclusions we require. We first state two lemmas which provide a mechanism for using a connected matching in one colour to build up a connected matching in the other colour.

\begin{lemma}[{\cite[Claim~6.3]{HaxellLuczakPengRodl}}] \label{lem:one-edge}
    Let $U=\{x,y,z,a,b,c,d\}\subseteq V(H)$ be a set of seven vertices in a~$2$-edge-coloured~$\eps$-complete~$3$-graph~$H=H_{\red}\sqcup H_{\blue}$. Suppose that
    \begin{enumerate} [label=(\roman{*})]
        \item $xyz$ is an edge of some red component $C_{\red}$,
        \item $ya,zb\in\partial C_{\blue}$ for some blue component $C_{\blue}$, and
        \item $\{xya,xzb,zbc,yad\}\subseteq E(H)$.
    \end{enumerate}
    Then at least one of the following holds.
    \begin{enumerate} [label=(\alph{*})]
        \item \label{item:two-red} $U$ contains two disjoint edges of $C_{\red}$.
        \item \label{item:one-cross-blue} $U$ contains an edge $e$ in $C_{\blue}$ such that $|e\cap\{a,b,c\}|=1$ and $|e\cap\{x,y,z\}|=2$.
        \item \label{item:two-blue} $U$ contains two disjoint edges of $C_{\blue}$.
    \end{enumerate}
\end{lemma}

\begin{proof}
    Suppose that~\ref{item:one-cross-blue} does not hold. Then we have~$xya,xzb\in C_{\red}$ and so~$ya,zb\in\partial C_{\red}$. Consider the edges~$zbc$ and~$yad$. If both are blue, then~$ya,zb\in\partial C_{\blue}$ implies~$zbc,yad\in C_{\blue}$, yielding~\ref{item:two-blue}. Otherwise, one of them, say~$zbc$, is red. Since~$zb\in\partial C_{\red}$, we have~$zbc\in C_{\red}$. Together with~$xya\in C_{\red}$, this yields~\ref{item:two-red}. The argument for~$yad$ is analogous.
\end{proof}

\begin{lemma}[{\cite[Claim~6.4]{HaxellLuczakPengRodl}}] \label{lem:two-edges}
    Let $U=\{u,v,w,x,y,z,a,b,c\}\subseteq V$ be a set of nine vertices in a~$2$-edge-coloured~$\eps$-complete~$3$-graph~$H=H_{\red}\sqcup H_{\blue}$. Suppose that
    \begin{enumerate} [label=(\roman{*})]
        \item $uvw$ and $xyz$ are edges of some red component $C_{\red}$,
        \item $ya,zb,vb,wc\in\partial C_{\blue}$ for some blue component $C_{\blue}$, and
        \item $\{xzb,vwc,uwc,yza,yua,xvb\}\subseteq E(H)$.
    \end{enumerate}
    Then at least one of the following holds.
    \begin{enumerate} [label=(\alph{*})]
        \item \label{item:three-red} $U$ contains three disjoint edges of $C_{\red}$.
        \item \label{item:pair-cross-blue} $U$ contains an edge $e$ in $C_{\blue}$ such that $|e\cap\{a,b,c\}|=1$ and either $|e\cap\{u,v,w\}|=2$ or $|e\cap\{x,y,z\}|=2$.
        \item \label{item:two-blue-nine} $U$ contains two disjoint edges of $C_{\blue}$.
    \end{enumerate}
\end{lemma}

\begin{proof}
    Suppose that~\ref{item:pair-cross-blue} does not hold. Then the edges~$xzb$, $vwc$, $uwc$ and~$yza$ are all red because~$ya,zb,wc\in\partial C_{\blue}$; in fact, they are all in~$C_{\red}$ because~$xyz,uvw\in C_{\red}$. Now consider the edges~$yua$ and~$xvb$. If both are blue, then~$ya,vb\in\partial C_{\blue}$ implies~$yua,xvb\in C_{\blue}$, yielding~\ref{item:two-blue-nine}. Otherwise, one of them, say~$yua$, is red. Since~$ya\in\partial C_{\red}$, we have~$yua\in C_{\red}$. Together with~$xzb,vwc\in C_{\red}$, this yields~\ref{item:three-red}. The argument for~$xvb$ is analogous.
\end{proof}

The next lemma is an easy generalisation of~\cite[Lemma~8.5]{HaxellLuczakPengRodl}. It relies heavily on the previous two lemmas and says that we may use a maximal matching $M$ in a red component to construct a suitably large matching $M'$ in a blue component. The key generalisation is that we remove the condition~$|P|\ge|M|+30\eps N$ by adjusting the size of~$M'$ we aim for.

\begin{lemma} \label{lem:mirror-matching}
    Let~$H=H_{\red}\sqcup H_{\blue}$ be a~$2$-edge-coloured~$\eps$-complete~$3$-graph on~$N$ vertices, $M$ be a largest matching in a red component $C_{\red}$ and $P$ be a set of vertices disjoint from $V(M)$. Let $q=\min(|M|,|P|-30\eps N)$. Suppose that there is a blue component $C_{\blue}$ such that for every $e\in M$ the bipartite induced subgraph $\partial C_{\blue}[e,P]$ contains $K_{2,|P|-9\eps N-1}$. Then there is a matching $M'$ in $C_{\blue}$ such that
    \begin{enumerate}[label=(\roman{*})]
        \item $|M'|\ge q$ and
        \item $|V(M')\cap P|\le q + 4\eps N$.
    \end{enumerate}
\end{lemma}

\begin{proof}
    Let~$M''$ be the largest matching in~$C_{\blue}$ such that
    \begin{equation} \label{eq:mirror-matching-interim}
        |V(M'')\cap P|\le |M''|\;\textrm{ and }\;|\{e\in M:e\cap V(M'')\ne\nth\}|\le|M''|\,.
    \end{equation}
    
    \begin{claim}
        $|M''|\ge q-4\eps N+2$.
    \end{claim}
    
    \begin{claimproof}
        Suppose that~$|M''|\le q-4\eps N+1$. We will show that there exist edges~$e_1,e_2\in M$ and vertices~$a,b,c\in P'':= P\sm V(M'')$ such that~$(e_1\cup e_2)\cap V(M'')=\nth$ and the set~$X=e_1\cup e_2\cup\{a,b,c\}$ satisfies the assumptions of Lemma~\ref{lem:two-edges}.

        By the second part of~\eqref{eq:mirror-matching-interim} and our supposed bound on~$|M''|$, there are at least~$4\eps N-1$ edges of~$M$ disjoint from~$V(M'')$. Let~$e_1=uvw\in M$ be any such edge. Since~$H$ is~$\eps$-complete, it follows from~\eqref{eq:complete-codegree} that for all~$x\in V(H)$ we have~$|N(x)|\ge(1-\eps)N+1$; this implies
        \[|V(H)\sm(N(u)\cap N(v)\cap N(w))| \le 3N-|N(u)|-|N(v)|-|N(w)|\le 3\eps N<4\eps N-2\,.\]
        Hence, we may pick an edge~$e_2=xyz\in M\sm\{e_1\}$ disjoint from~$V(M'')$ such that every pair of vertices in~$\binom{e_1\cup e_2}{2}$ is active. 

        By the first part of~\eqref{eq:mirror-matching-interim} and our bounds on~$|P|$ and~$|M''|$, we have
        \[|P''|=|P|-|P\cap V(M'')| \ge q + 30\eps N - (q - 4\eps N+1) = 34\eps N-1\,.\]
        Both~$\partial C_{\blue}[e_1,P]$ and~$\partial C_{\blue}[e_2,P]$ contain a copy~$K_{2,|P|-9\eps N-1}$ by assumption, and the number of vertices of~$P''$ not in both copies of~$K_{2,|P|-9\eps N-1}$ is at most~$18\eps N+2$. Every pair of vertices in~$\binom{e_1\cup e_2}{2}$ is active, so by~\eqref{eq:complete-codegree} we have
        \[\Bigg|P''\sm\Bigg(\bigcap_{ij}N(i,j)\Bigg)\Bigg| \le 15\eps N\,,\]
        where the intersection is taken over all pairs~$ij\in\binom{e_1\cup e_2}{2}$. Since~$34\eps N- 18\eps N-15\eps N - 3 \ge 3$, we may choose distinct vertices~$a,b,c\in P''$ so that
        \begin{enumerate} [label=(\roman{*})]
            \item $\left\{e\in\binom{e_1\cup e_2\cup\{a,b,c\}}{3}:|e\cap(e_1\cup e_2)|=2\right\}\subseteq E(H)$ and
            \item $K_{2,3}\subseteq\partial C_{\blue}[e_i,\{a,b,c\}]$ for~$i=1,2$.
        \end{enumerate}

        Hence, we may apply Lemma~\ref{lem:two-edges} to the set~$X=e_1\cup e_2\cup\{a,b,c\}$. Since we may enlarge~$M$ in~$C_{\red}$ if~\ref{item:three-red} of Lemma~\ref{lem:two-edges} occurs, the maximality of~$M$ implies that~\ref{item:pair-cross-blue} or~\ref{item:two-blue-nine} of Lemma~\ref{lem:two-edges} occurs. But now we may enlarge~$M''$ in~$C_{\blue}$ with the conditions~\eqref{eq:mirror-matching-interim} preserved, yielding a contradiction with the choice of~$M''$.
    \end{claimproof}

    We obtain~$M'$ from~$M''$ as follows. If~$|M''|\ge q$, we simply let~$M'=M''$. Otherwise, we perform the following iterative procedure which enlarges~$M''$ by two in each iteration until its size reaches~$q$. Write~$M'$ for the current matching; in each step we maintain~$M''\subseteq M'$, add at most four vertices of~$P$ to~$V(M')$ and maintain~$|\{e\in M:e\cap V(M')\ne\nth\}|\le|M'|$. Since there are~$\left\lceil\frac{q-|M''|}{2}\right\rceil$ steps, for the final~$M'$ we have
    \[|P\cap V(M')| \le |M''|+ 2(q-|M''|+1) \le q + 4\eps N\,.\]
    Now we describe a single step of the procedure. Take~$e=xyz\in M$ such that~$e\cap V(M')=\nth$.

    Fix a copy of~$K_{2,|P|-9\eps N-1}$ in~$\partial C_{\blue}[e,P]$ (which exists by assumption) and let~$P_0$ denote the set of~$9\eps N +1$ vertices in~$P$ not in this copy of~$K_{2,|P|-9\eps N-1}$. Set~$P'=P\sm(V(M')\cup P_0)$. By our bounds on~$|P|$ and~$|P\cap V(M')|$, we have
    \[|P'|\ge q+30\eps N - (q+4\eps N) - |P_0| \ge 17\eps N-1\,.\]
    Set~$N_1=P'\cap N(x,y)\cap N(x,z)\cap N(y,z)$; by~\eqref{eq:complete-codegree} we have~$|N_1|\ge|P'|-3\eps N\ge14\eps N-1$, so we may pick~$a\in N_1$. Set~$N_2=N_1\cap N(a,x)\cap N(a,y)\cap N(a,z)$; by~\eqref{eq:complete-codegree} we have~$|N_2|\ge|N_1|-3\eps N\ge11\eps N-1$, so we may pick~$b\in N_2$. Set~$N_3=N_2\cap N(b,x)\cap N(b,y)\cap N(b,z)$; by~\eqref{eq:complete-codegree} we have~$|N_3|\ge|N_2|-3\eps N\ge8\eps N-1$, so we may pick~$c\in N_3$. Finally, set~$N_4=N_3\cap N(c,x)\cap N(c,y)\cap N(c,z)$; by~\eqref{eq:complete-codegree} we have~$|N_4|\ge|N_3|-3\eps N\ge5\eps N-1$, so we may pick~$d\in N_4$. Hence, we may pick~$a,b,c,d\in P'$ so that
    \begin{enumerate} [label=(\roman{*})]
        \item $K_{2,3}\subseteq\partial C_{\blue}[e,\{a,b,c\}]$ and
        \item $\left\{f\in\binom{e\cup\{a,b,c,d\}}{3}:|f\cap e|\in\{1,2\}\right\}\subseteq E(H)$.
    \end{enumerate}
    Now we apply Lemma~\ref{lem:one-edge} to the set~$X=e\cup\{a,b,c,d\}$. If~\ref{item:two-red} of Lemma~\ref{lem:one-edge} occurs, we may enlarge~$M$ in~$C_{\red}$; if~\ref{item:one-cross-blue} of Lemma~\ref{lem:one-edge} occurs, we may enlarge~$M''$ in~$C_{\blue}$ with the conditions~\eqref{eq:mirror-matching-interim} preserved. Now the maximality of~$M$ and~$M''$ implies that~\ref{item:two-blue} of Lemma~\ref{lem:one-edge} occurs, so by adding the two edges it guarantees, we may enlarge~$M'$ by two and use just one edge of~$M$ and four vertices of~$P$.
\end{proof}

We are now ready to prove Lemma~\ref{lem:matching-div-div}.

\begin{proof}[Proof of Lemma~\ref{lem:matching-div-div}]
    By Lemmas~\ref{lem:vertex-shadow-monoxcomp} and~\ref{lem:dominant-monoxcomp}, we may assume Setup~\ref{setup} with~$N=3s+t+\eta s$ from now onwards and at least one of~$S$ and~$A$ exists. Let~$M$ be a largest matching in~$H$ among all matchings contained in~$S$ or~$A$; without loss of generality, we may assume that~$M\ne\nth$. The arguments are largely identical, so we shall assume~$M\subseteq S$ and explain as needed the adjustments required for~$M\subseteq A$ (with the colours interchanged). Suppose that~$1 \le m := |M| < s$; for~$M\subseteq A$ suppose that~$1 \le m < t$. Set~$R:=V'_{\red}\sm V(M)$ and~$B:=V'_{\blue}\sm V(M)$. Note that~$R\cap B=\nth$ and by Lemma~\ref{lem:setup-properties}\ref{item:colour-bipartition-near-complete} we have
    \begin{equation} \label{eq:remainder-size}
        |R\cup B| = N - 3m - |V(H)\sm V'| \ge 3s + t + \eta s - 3m - 8\eps N\,.
    \end{equation}

    \begin{claim} \label{claim:same-blue-monoxcomp}
        If~$|R|>2\eps N$, then all edges~$xyz\in H[R\cup B]$ with~$x\in R$ and~$xy\in\partial S$ belong to the same blue component~$C_{\blue}$. Furthermore, if also~$|B|>2\eps N$, then~$C_{\blue}=A$.
    \end{claim}

    \begin{claimproof}
        Let~$G$ be the graph on~$R\cup B$ such that~$xy\in E(G)$ if and only if~$x\in R$ and~$xy\in\partial S$, and obtain~$G'$ from~$G$ by deleting all isolated vertices. First, note that every edge~$xyz\in H[R\cup B]$ with~$xy\in E(G')$ is blue. Indeed, any such edge which is red would be a red edge in~$S$ disjoint from~$M$, contradicting the maximality of~$M$. For~$xy\in E(G')$ let~$C_{xy}$ be the unique blue component containing every edge~$xyz$ with~$z\in N_H(x,y)\cap(R\cup B)$; this is well-defined by definition, the observation above and~\eqref{eq:complete-codegree}. Now for adjacent edges~$xy$ and~$xz$ of~$G'$, observe that since~$xy,xz\in\partial S$ and~$|R\cup B|\ge2\eps N$, by~\eqref{eq:complete-codegree} there exists~$w\in R\cup B$ such that~$xyw,xzw\in E(H)$. In particular, $C_{xy}=C_{xz}$. But now by Lemma~\ref{lem:setup-properties}\ref{item:dominant-monoxcomp-shadows-nbrhood-sizes} the subgraph~$G'[R]$ has minimum degree at least~$|R|-\eps N\ge|R|/2$, so it is connected; it follows by definition that~$G'$ itself is also connected. Hence, there is a blue component~$C_{\blue}$ such that for all~$xy\in E(G)$ we have~$C_{xy}=C_{\blue}$.

        To prove the second part, note that since~$|R|,|B|> 2\eps N$, by Lemma~\ref{lem:setup-properties}\ref{item:dominant-monoxcomp-shadows-nbrhood-sizes} more than~$\frac{|R||B|}{2}$ edges of~$\partial S$ have one end vertex in~$R$ and the other in~$B$, and the same holds for the edges of~$\partial A$. Hence, there exist~$x\in R$ and~$y\in B$ such that~$xy\in\partial S\cap \partial A$, so we have~$C_{\blue}=A$.
    \end{claimproof}
    
    For the rest of the proof we distinguish three cases.

    \textbf{Case 1:}~$|B|\le5\eps N$. Let~$C_{\blue}$ denote the blue component returned by Claim~\ref{claim:same-blue-monoxcomp}. 

    \begin{claim} \label{claim:dense-shadow-small-B}
        For each edge~$e\in M$ the bipartite induced subgraph~$\partial C_{\blue}[e,R]$ contains a copy of~$K_{2,|R|-3\eps N-1}$.
    \end{claim}

    \begin{claimproof}
        Let~$e\in M$. Write~$R_e$ for the set of vertices~$a\in R$ such that~$xa,ya,za\in\partial H$; by~\eqref{eq:complete-codegree} we have~$|R_e|\ge|R|-3\eps N$. Suppose that~$\partial C_{\blue}[e,R_e]$ does not contain a copy of~$K_{2,|R_e|-1}$. Then there are distinct vertices~$a,b\in R_e$ such that~$ya,zb\notin\partial C_{\blue}$. Since~$|R|\ge2\eps N+5$, by~\eqref{eq:complete-codegree} and Lemma~\ref{lem:setup-properties}\ref{item:dominant-monoxcomp-shadows-nbrhood-sizes} we may pick vertices~$c,d,u\in R\sm\{a,b\}$ such that~$ac,bd\in\partial S$ and~$yac,zbd,uac,ubd\in E(H)$. By Claim~\ref{claim:same-blue-monoxcomp} we have~$uac,ubd\in C_{\blue}$, so~$ac,bd\in\partial C_{\blue}$. Since~$ya,zb\notin\partial C_{\blue}$ and~$ac,bd\in\partial S$, the edges~$yac$ and~$zbd$ must be red edges in~$S$. But then~$(M\sm\{xyz\})\cup\{yac,zbd\}$ is a larger matching in~$S$ than $M$, which is a contradiction.
    \end{claimproof}

    Since Claim~\ref{claim:dense-shadow-small-B} holds, we may apply Lemma~\ref{lem:mirror-matching} with~$C_{\red}=S$, $P=R$ and~$C_{\blue}$ to obtain a matching~$M'$ in~$C_{\blue}$ of size~$|M'|=:m'\ge q:=\min(m,|R|-30\eps N)$ with~$|R\cap V(M')|\le q+4\eps N$. Note that we are done if~$m'\ge s$. Otherwise, by~\eqref{eq:remainder-size} and our choice of constants we have
    \begin{equation} \label{eq:big-R-small-B-size}
        |R| \ge 3s + t +\eta s - 3m - 13\eps N \ge
        \begin{cases}
            m + 30\eps N &\textrm{ if }\;m<t \\
            t + 30\eps N &\textrm{ if }\;t\le m<s\,.
        \end{cases}
    \end{equation}
    Note that if~$t\le m< s$ then~$m'\ge q\ge t$ and we are done; for~$M\subseteq A$ this case does not happen. Hence, it remains to consider the case~$m<t$. In this case, by~\eqref{eq:big-R-small-B-size} we have
    \begin{equation*}
        |R\sm V(M')| \ge 3s + t +\eta s - 4m - 17\eps N \ge 3(s-m') + 3\eps N\,.
    \end{equation*}
    Now Lemma~\ref{lem:setup-properties}\ref{item:dominant-monoxcomp-shadows-nbrhood-sizes} and~\eqref{eq:complete-codegree} allow us to greedily find the required disjoint edges~$xyz\in H[R]$ with~$xy\in\partial S$. In the current case~$M\subseteq S$ we only need~$t-m'$ such edges; in the case~$M\subseteq A$, the same estimate gives the full~$s-m'$ edges. Now by Claim~\ref{claim:same-blue-monoxcomp} all these edges are in~$C_{\blue}$, so we may add them to~$M'$ to obtain a matching of size~$s$ in~$C_{\blue}$.

    \textbf{Case 2:}~$|R|\le5\eps N$. By our assumptions and~\eqref{eq:remainder-size}, we have~$B\ne\nth$ and so~$A$ exists.

    \begin{claim} \label{claim:dense-shadow-small-R}
        For each edge~$e\in M$ the bipartite induced subgraph~$\partial A[e,B]$ contains a copy of~$K_{2,|B|-9\eps N-1}$.
    \end{claim}

    \begin{claimproof}
        Fix an edge~$e=xyz\in M$. Write~$B_e$ for the set of vertices~$a\in B$ such that~$xa,ya,za\in\partial H$; by~\eqref{eq:complete-codegree} we have~$|B_e|\ge|B|-3\eps N$. Call a vertex~$a\in B_e$ \emph{friendly} to~$x$ if~$xa\in \partial S\cup \partial A$, and let~$B_x$ be the subset of~$B_e$ containing all unfriendly vertices to~$x$.

        \begin{claim}
            $|B_x|<2\eps N$.
        \end{claim}

        \begin{claimproof}
        Suppose that~$|B_x| \ge 2\eps N$. Since~$S$ exists, recall that~$|V'_{\red}|\ge4\eps N$. Let~$G_S$ and~$G_A$ denote the bipartite induced subgraphs of~$\partial S$ and~$\partial A$, respectively, on vertex set~$B_x\cup V'_{\red}$. Taking subsets if necessary, assume for simplicity that~$|B_x|=2\eps N$ and~$|V'_{\red}|=4\eps N$. By Lemma~\ref{lem:setup-properties}\ref{item:dominant-monoxcomp-shadows-nbrhood-sizes} we have~$|E(G_S)|\ge4\eps^2N^2$ and~$|E(G_A)|\ge6\eps^2N^2$, so~$|E(G_S)\cap E(G_A)|\ge2\eps^2N^2$. Hence, we may pick~$a\in B_x$ with degree at least~$\eps N$ in~$G_S\cap G_A$. Then, by the definition of~$B_e$ and~\eqref{eq:complete-codegree}, there is a vertex~$b\in V'_{\red}$ such that~$xab\in E(H)$ and~$ab\in E(G_S)\cap E(G_A)$. But this implies that~$xa\in \partial S\cup \partial A$, no matter how~$xab$ is coloured, thereby giving a contradiction.
        \end{claimproof}

        Set~$B'_e:=B_e\sm(B_x\cup B_y\cup B_z)$, and suppose that~$\partial A[e,B'_e]$ contains no copy of~$K_{2,|B'_e|-1}$. Then there are distinct vertices~$a,b\in B'_e$ such that~$ya,zb\notin\partial A$; in particular, $ya,zb\in\partial S$. Since~$|B|\ge2\eps N+4$, by~\eqref{eq:complete-codegree} and Lemma~\ref{lem:setup-properties}\ref{item:dominant-monoxcomp-shadows-nbrhood-sizes} we may pick vertices~$c,d\in B\sm\{a,b\}$ such that~$ac,bd\in\partial A$ and~$yac,zbd\in E(H)$. Hence, the edges~$yac$ and~$zbd$ must be red and in~$S$. But then~$(M\sm\{xyz\})\cup\{yac,zbd\}$ is a larger matching in~$S$ than $M$, which is a contradiction.
    \end{claimproof}

    Since Claim~\ref{claim:dense-shadow-small-R} holds, we may apply Lemma~\ref{lem:mirror-matching} with~$C_{\red}=S$, $P=B$ and~$C_{\blue}=A$ to obtain a matching~$M'$ in~$A$ of size~$|M'|=:m'\ge q:=\min(m,|B|-30\eps N)$ with~$|B\cap V(M')|\le q+4\eps N$. Note that we are done if~$m'\ge s$. Otherwise, by~\eqref{eq:remainder-size} and our choice of constants we have
    \begin{equation} \label{eq:big-B-small-R-size}
        |B| \ge 3s + t +\eta s - 3m - 13\eps N \ge
        \begin{cases}
            m + 30\eps N &\textrm{ if }\;m<t \\
            t + 30\eps N &\textrm{ if }\;t\le m<s\,.
        \end{cases}
    \end{equation}
    Note that if~$t\le m< s$ then~$m'\ge q\ge t$ and we are done; for~$M\subseteq A$ this case does not happen. Hence, it remains to consider the case~$m<t$. In this case, the maximality of~$M$ implies~$m'=m$.
    
    We claim that~$R':=V'_{\red}\sm V(M')=\nth$. Indeed, suppose that there is~$d\in R'$. By~\eqref{eq:big-B-small-R-size} we have
    \begin{equation*}
        |B\sm V(M')| \ge 3s + t +\eta s - 4m - 17\eps N \ge 2\eps N + 3\,,
    \end{equation*}
    so by~\eqref{eq:complete-codegree} and Lemma~\ref{lem:setup-properties}\ref{item:dominant-monoxcomp-shadows-nbrhood-sizes} we may choose vertices~$a,b,c\in B\sm V(M')$ such that~$ad\in\partial S$, $ab\in\partial A$ and~$abd,abc\in E(H)$. 
    Consider the edges~$abd$ and~$abc$. If both are red, then~$ad\in\partial S$ implies~$abd,abc\in S$ and~$M\cup\{abc\}$ is a larger matching in~$S$, which is a contradiction. Hence, one of them, say~$abc$, is blue. But then~$ab\in\partial A$ implies that~$abc\in A$ and so~$M'\cup\{abc\}$ is a larger matching in~$A$, which yields a contradiction. Hence, $R'=\nth$ and we are in Case 1 with the colours red and blue interchanged and the matching~$M$ replaced by~$M'$.

    \textbf{Case 3:}~$|R|,|B|\ge5\eps N$.

    \begin{claim} \label{claim:dense-shadow-R-B}
        For every edge~$e\in M$ the bipartite induced subgraph~$\partial A[e,R\cup B]$ contains a copy of~$K_{2,|R\cup B|-9\eps N-1}$.
    \end{claim}

    \begin{claimproof}
        Fix an edge~$e=xyz\in M$. Write~$P_e$ for the set of vertices~$a\in R\cup B$ such that~$xa,ya,za\in\partial H$; by~\eqref{eq:complete-codegree} we have~$|P_e|\ge|R\cup B|-3\eps N$. Call a vertex~$a\in P_e \cap B$ \emph{friendly} to~$x$ if~$xa\in \partial S\cup \partial A$, and let~$B_x$ be the subset of~$P_e\cap B$ containing all unfriendly vertices to~$x$.

        \begin{claim}
            $|B_x|<2\eps N$.
        \end{claim}

        \begin{claimproof}
        Suppose that~$|B_x| \ge 2\eps N$. Since~$S$ exists, recall that~$|V'_{\red}|\ge4\eps N$. Let~$G_S$ and~$G_A$ denote the bipartite induced subgraphs of~$\partial S$ and~$\partial A$, respectively, on vertex set~$B_x\cup V'_{\red}$. Taking subsets if necessary, assume for simplicity that~$|B_x|=2\eps N$ and~$|V'_{\red}|=4\eps N$. By Lemma~\ref{lem:setup-properties}\ref{item:dominant-monoxcomp-shadows-nbrhood-sizes} we have~$|E(G_S)|\ge4\eps^2N^2$ and~$|E(G_A)|\ge6\eps^2N^2$, so~$|E(G_S)\cap E(G_A)|\ge2\eps^2N^2$. Hence, we may pick~$a\in B_x$ with degree at least~$\eps N$ in~$G_S\cap G_A$. Then, by the definition of~$P_e$ and~\eqref{eq:complete-codegree}, there is a vertex~$b\in V'_{\red}$ such that~$xab\in E(H)$ and~$ab\in E(G_S)\cap E(G_A)$. But this implies that~$xa\in \partial S\cup \partial A$, no matter how~$xab$ is coloured, thereby giving a contradiction.
        \end{claimproof}

        Set~$P'_e:=P_e\sm(B_x\cup B_y\cup B_z)$, and suppose that~$\partial A[e,P'_e]$ does not contain a copy of~$K_{2,|P'_e|-1}$. Then there are distinct vertices~$a,b\in P'_e$ such that~$ya,zb\notin\partial A$. The rest of the proof combines arguments from the proofs of Claims~\ref{claim:dense-shadow-small-B} and~\ref{claim:dense-shadow-small-R}, splitting into cases depending on whether each of~$a$ and~$b$ is in~$R$ or~$B$. Set~$d_a:=y$ and~$d_b:=z$. For~$w\in\{a,b\}\cap R$, we may pick distinct vertices~$c_w,u_w\in R\sm\{a,b\}$ such that~$wc_w\in\partial S$ and~$d_wwc_w,wc_wu_w\in E(H)$. By Claim~\ref{claim:same-blue-monoxcomp} we have~$wc_wu_w\in A$, so~$wc_w\in\partial A$. Since~$d_ww\notin\partial A$ and~$wc_w\in\partial S$, the edge~$d_wwc_w$ must be red and in~$S$. For~$w\in\{a,b\}\cap B$, we may pick a distinct vertex~$c_w\in B\sm\{a,b\}$ such that~$wc_w\in\partial A$ and~$d_wwc_w\in E(H)$. Hence, the edge~$d_wwc_w$ must be red and in~$S$. Since~$|R|,|B|\ge2\eps N+6$, by~\eqref{eq:complete-codegree} and Lemma~\ref{lem:setup-properties}\ref{item:dominant-monoxcomp-shadows-nbrhood-sizes} the aforementioned choices are possible. But then~$(M\sm\{xyz\})\cup\{yac_a,zbc_b\}$ is a larger matching in~$S$ than $M$, which is a contradiction.
    \end{claimproof}

    Fix~$a\in R$ and~$b\in B$ such that~$ab\in\partial H$; such a pair exists because~$|R|,|B|\ge5\eps N$. Since Claim~\ref{claim:dense-shadow-R-B} holds, we may apply Lemma~\ref{lem:mirror-matching} with~$C_{\red}=S$, $P=(R\cup B)\sm\{a,b\}$ and~$C_{\blue}=A$ to obtain a matching~$M'$ in~$A$ of size~$|M'|=:m'\ge q:=\min(m,|P|-30\eps N)$ with~$|P\cap V(M')|\le q+4\eps N$. Note that we are done if~$m'\ge s$. Otherwise, by~\eqref{eq:remainder-size} and our choice of constants we have
    \begin{equation} \label{eq:R-B-size}
        |P| \ge 3s + t +\eta s - 3m - 14\eps N -2 \ge
        \begin{cases}
            m + 30\eps N &\textrm{ if }\;m<t \\
            t + 30\eps N &\textrm{ if }\;t\le m<s\,.
        \end{cases}
    \end{equation}
    Note that if~$t\le m< s$ then~$m'\ge q\ge t$ and we are done; for~$M\subseteq A$ this case does not happen. Hence, it remains to consider the case~$m<t$. In this case, note that~$m'=m$ by the maximality of~$M$. By~\eqref{eq:R-B-size} we have
    \begin{equation*}
        |P\sm V(M')| \ge 3s + t +\eta s - 4m - 18\eps N - 2 \ge 3\eps N + 1\,,
    \end{equation*}
    so by~\eqref{eq:complete-codegree} and Lemma~\ref{lem:setup-properties}\ref{item:dominant-monoxcomp-shadows-nbrhood-sizes} we may choose a vertex~$c\in P\sm V(M')$ so that~$ac\in\partial S$, $bc\in\partial A$ and~$abc\in E(H)$. Observe that~$abc$ is an edge disjoint from both~$V(M)$ and~$V(M')$ which is in~$S$ if it is red and in~$A$ if it is blue. Hence, either~$M\cup\{abc\}$ is a matching in~$S$ which is larger than~$M$, or~$M'\cup\{abc\}$ is a matching in~$A$ which is larger than~$M$, yielding a contradiction.
\end{proof}

\section{Proving lemmas} \label{sec:lemma-proofs-matching-xdivs}

In this appendix we provide the proofs of several lemmas from Section~\ref{sec:matching-xdivs}. We first state a technical lemma which is useful in these proofs.

\begin{lemma}[{\cite[Fact 7.9]{HaxellLuczakPengRodl}}] \label{lem:parts-allocations}
    Given numbers~$a_1\ge\dots\ge a_k\ge0$, let~$T:=\sum_{i\in[k]}a_i$.
    \begin{enumerate}[label=(\roman{*})]
        \item \label{item:two-parts} Let~$k\ge2$ and~$d\ge\frac{2}{3}T$. If~$a_1\le d$, then there exists~$\ell_0\in[k-1]$ such that~$T-d \le \sum_{i\in[\ell_0]}a_i \le d$.
        \item \label{item:three-parts} Let~$k\ge3$ and~$T\ge5r$. If~$a_1+a_2\le T-2r$, then there exist~$1\le\ell_1<\ell_2\le k-1$ such that~$\sum_{i\in[\ell_1]}a_i \ge r$, $\sum_{i=\ell_1+1}^{\ell_2}a_i \ge r$ and~$\sum_{i=\ell_2+1}^{k}a_i \ge r$.
    \end{enumerate}
\end{lemma}

\begin{proof}[Proof of Lemma~\ref{lem:vertex-shadow-monoxcomp}]
    Fix~$x\in V(H)$, and let~$H(x)$ denote the graph with vertex set~$V(H)\sm\{x\}$ and edge set~$\{yz:xyz\in E(H)\}$; note that~$H(x)$ has minimum degree at least~$(1-\eps)N\ge\frac{3N}{4}$ and the~$2$-edge-colouring~$H=H_{\red}\sqcup H_{\blue}$ induces a~$2$-edge-colouring~$H(x)=H_{\red}(x)\sqcup H_{\blue}(x)$.
    
    \begin{claim} \label{claim:link-large-monoxcomp}
        $H(x)$ contains a monochromatic component with at least~$(1-\eps)N$ vertices.
    \end{claim}

    \begin{claimproof}
            Suppose that every component of~$H_{\red}(x)$ has fewer than~$(1-\eps)N$ vertices. Then by Lemma~\ref{lem:parts-allocations}\ref{item:two-parts} there is a partition~$V(H(x))=V_1\sqcup V_2$ such that~$\eps N\le|V_1|\le|V_2|< (1-\eps)N$ and~$E_{H_{\red}(x)}(V_1,V_2)=\nth$. Now the minimum degree condition of~$H(x)$ implies that every vertex in~$V_2$ has a neighbour in~$V_1$ and every vertex in~$V_1$ has more than~$\frac{|V_2|}{2}$ neighbours in~$V_2$, so the graph~$H_{\blue}(x)$ is connected and there is a blue component on all~$N-1\ge(1-\eps)N$ vertices.
    \end{claimproof}

    Claim~\ref{claim:link-large-monoxcomp} returns a monochromatic component in~$H(x)$ with at least~$(1-\eps)N$ vertices. This naturally corresponds to a monochromatic component~$C_x$ in~$H$ which satisfies~\eqref{eq:vertex-shadow-monoxcomp}.
\end{proof}

\begin{proof}[Proof of Lemma~\ref{lem:dominant-monoxcomp}]
    Let~$G$ be the graph on~$V_{\red}$ with edge set~$\left\{xy\in\binom{V_{\red}}{2}:xy\in\partial C_x\cap\partial C_y\right\}$.

    \begin{claim} \label{claim:many-same-monoxcomp}
        $G$ contains a component with at least~$|V_{\red}|-2\eps N$ vertices.
    \end{claim}

    \begin{claimproof}
        Suppose that every component of~$G$ has fewer than~$|V_{\red}|-2\eps N$ vertices. Then by Lemma~\ref{lem:parts-allocations}\ref{item:two-parts} there is a partition~$V_{\red}=V_1\sqcup V_2$ such that~$2\eps N \le |V_1| \le |V_2| < |V_{\red}|-2\eps N$ and~$E_{G}(V_1,V_2)=\nth$. Hence, $G$ has at least~$2\eps N(|V_{\red}|-2\eps N) > |V_{\red}|\eps N$ nonedges. But~\eqref{eq:vertex-shadow-monoxcomp} implies that~$G$ has at most~$|V_{\red}|\eps N$ nonedges, which is a contradiction.
    \end{claimproof}

    Claim~\ref{claim:many-same-monoxcomp} guarantees that~$G$ has a component~$C$ with at least~$|V_{\red}|-2\eps N$ vertices. Since by~\eqref{eq:components-disjoint-shadows}~$xy\in\partial C_x\cap\partial C_y$ implies~$C_x=C_y$, there is a red component~$C_{\red}$ such that~$C_x=C_{\red}$ for all~$x\in V(C)$. This completes the proof.
\end{proof}

\begin{proof}[Proof of Lemma~\ref{lem:dominant-pair-small-red-set-blue-C4}]
    We first establish a version of~\eqref{eq:blue-monoxcomps-dominant-pair} involving the azure component~$A$ whenever it exists.

    \begin{claim} \label{claim:blue-monoxcomps-dominant-pair-azure}
        If~$V'_{\blue}\ne\nth$, then there is a blue component~$C_{\blue}$ such that
        \begin{equation} \label{eq:blue-monoxcomps-dominant-pair-azure}
            \left|\binom{V'_{\red}}{2}\sm\left(\partial A\cup\partial C_{\blue}\right)\right| < 4\eps^{1/6}N^2\,.
        \end{equation}
    \end{claim}

    \begin{claimproof}
        If~$\left|\partial C^{(i)}_{\blue}\cap\binom{V'_{\red}}{2}\right| \le 3\eps^{1/6}N^2$ holds for some~$i\in\{1,2\}$, then by taking~$C_{\blue}=C^{(3-i)}_{\blue}$ and applying~\eqref{eq:blue-monoxcomps-dominant-pair} we have
        \begin{equation*}
            \left|\binom{V'_{\red}}{2}\sm\left(\partial A\cup\partial C_{\blue}\right)\right| < 2\eps^{1/3}N^2 + \left|\partial C^{(i)}_{\blue}\cap\binom{V'_{\red}}{2}\right| \le 4\eps^{1/6}N^2\,.
        \end{equation*}
        Hence, suppose that~$\left|\partial C^{(i)}_{\blue}\cap\binom{V'_{\red}}{2}\right| > 3\eps^{1/6}N^2$ for~$i\in\{1,2\}$. We prove the following subclaim.
        
        \begin{claim}
            There are vertices~$a\in V'_{\blue}$ and~$u,v,w\in V'_{\red}$ such that~$H[\{a,u,v,w\}]$ is a copy of~$C^{(3)}_4$, $au,av,aw\in\partial A$, and setting~$P:=\{uv,uw,vw\}$, we have~$P\cap\partial S\ne\nth$, $P\cap\partial C^{(1)}_{\blue}\ne\nth$, $P\cap\partial C^{(2)}_{\blue}\ne\nth$ and~$P\subseteq\partial C^{(1)}_{\blue}\cup\partial C^{(2)}_{\blue}$.
        \end{claim}

        \begin{claimproof}
            Pick~$a\in V'_{\blue}$ arbitrarily. Every graph~$G$ contains a vertex with degree at least~$\frac{2|E(G)|}{|V(G)|}$, so by~\ref{item:few-blue-vertices}, Lemma~\ref{lem:setup-properties}\ref{item:colour-bipartition-near-complete} and~\eqref{eq:blue-monoxcomps-dominant-pair} there are at least~$2\eps N$ vertices~$u\in V'_{\red}$ such that
            \begin{equation} \label{eq:shadow-nbr-dominant-pair}
                |\partial C^{(1)}_{\blue}(u)\cap V'_{\red}| + |\partial C^{(2)}_{\blue}(u)\cap V'_{\red}| > |V'_{\red}| - \eps^{1/4}N\,.
            \end{equation}
            Hence, by Lemma~\ref{lem:setup-properties}\ref{item:dominant-monoxcomp-shadows-nbrhood-sizes} there is a vertex~$u\in V'_{\red}$ satisfying~\eqref{eq:shadow-nbr-dominant-pair} such that~$au\in\partial A$.

            Since~$au\in\partial A\subseteq\partial H$, by~\eqref{eq:complete-codegree} there are at most~$4\eps N^2$ pairs~$vw$ such that~$H[\{a,u,v,w\}]$ is not a copy of~$C^{(3)}_4$. Because~$a\in V'_{\blue}$ and by Lemma~\ref{lem:setup-properties}\ref{item:dominant-monoxcomp-shadows-nbrhood-sizes}, there are at most~$2\eps N^2$ pairs~$vw$ such that~$av\notin\partial A$ or~$aw\notin\partial A$; analogously, because~$u\in V'_{\red}$ and by Lemma~\ref{lem:setup-properties}\ref{item:dominant-monoxcomp-shadows-nbrhood-sizes}, there are at most~$2\eps N^2$ pairs~$vw$ such that~$uv\notin\partial S$ or~$uw\notin\partial S$. By~\eqref{eq:blue-monoxcomps-dominant-pair} there are at most~$2\eps^{1/3}N^2$ pairs~$vw\in\binom{V'_{\red}}{2}$ such that~$vw\notin\partial C^{(1)}_{\blue}\cup\partial C^{(2)}_{\blue}$. Writing~$Y\subseteq\binom{V'_{\red}}{2}$ for the set of all such exceptional pairs, we have~$|Y|\le3\eps^{1/3}N^2$.

            If there are vertices~$v\in \partial C^{(1)}_{\blue}(u)\cap V'_{\red}$ and~$w\in \partial C^{(2)}_{\blue}(u)\cap V'_{\red}$ such that~$vw\notin Y$, then there are vertices~$a\in V'_{\blue}$ and~$u,v,w\in V'_{\red}$ with the desired properties. Now suppose that no such vertices~$v$ and~$w$ exist; in particular, we have~$\prod_{i=1}^2|\partial C^{(i)}_{\blue}(u)\cap V'_{\red}| \le |Y| \le 3\eps^{1/3}N^2$. Then for some~$i\in\{1,2\}$ we have~$|\partial C^{(i)}_{\blue}(u)\cap V'_{\red}| \le \sqrt{3}\eps^{1/6}N$; by~\eqref{eq:shadow-nbr-dominant-pair} we have~$|\partial C^{(3-i)}_{\blue}(u)\cap V'_{\red}| \ge |V'_{\red}| - 2\eps^{1/6}N$. It follows that at most~$2\eps^{1/6}N|V'_{\red}|\le 2\eps^{1/6}N^2$ pairs of~$\partial C^{(i)}_{\blue}\cap\binom{V'_{\red}}{2}$ are not contained in~$\binom{\partial C^{(3-i)}_{\blue}(u)\cap V'_{\red}}{2}$. Now since
            \[\left|\partial C^{(i)}_{\blue}\cap\binom{V'_{\red}}{2}\right|>3\eps^{1/6}N^2 \ge 2\eps^{1/6}N^2 + |Y |\,,\]
            there are~$v,w\in\partial C^{(3-i)}_{\blue}(u)\cap V'_{\red}$ such that~$vw\in\partial C^{(i)}_{\blue}\sm Y$. Hence, there are vertices~$a\in V'_{\blue}$ and~$u,v,w\in V'_{\red}$ with the desired properties.
        \end{claimproof}

        Since~$P\cap\partial C^{(1)}_{\blue}\ne\nth$ and~$P\cap\partial C^{(2)}_{\blue}\ne\nth$, the edge~$uvw$ must be red. Now~$P\cap\partial S\ne\nth$ and~$C^{(3)}_4\nsubseteq S$ by~\ref{item:scarlet-no-C4}, so one of the edges~$uva,uwa,vwa$ must be blue. Since~$au,av,aw\in\partial A$ and~$P\subseteq\partial C^{(1)}_{\blue}\cup\partial C^{(2)}_{\blue}$, there exists~$i\in[2]$ such that~$C^{(i)}_{\blue}=A$. The proof is completed by setting~$C_{\blue} = C^{(3-i)}_{\blue}$ and recalling~\eqref{eq:blue-monoxcomps-dominant-pair}.
    \end{claimproof}

    By Claim~\ref{claim:blue-monoxcomps-dominant-pair-azure}, from now onwards we may assume that
    \begin{equation} \label{eq:blue-monoxcomps-dominant-pair-adjusted}
        \left|\binom{V'_{\red}}{2}\sm\left(\partial C^{(1)}_{\blue}\cup\partial C^{(2)}_{\blue}\right)\right| < 4\eps^{1/6}N^2\,,
    \end{equation}
    and that~$C^{(1)}_{\blue}=A$ whenever~$A$ exists. To complete the proof, fix~$X\subseteq V'_{\red}$ with~$|X|\ge25\eps^{1/12}N$. By~\eqref{eq:blue-monoxcomps-dominant-pair-adjusted} there is a subset~$X'\subseteq X$ with~$|X\sm X'|\le8\eps^{1/12}N$ such that for each~$u\in X'$ there are at most~$\eps^{1/12}N$ vertices~$v\in X'$ such that~$uv\notin\partial C^{(1)}_{\blue}\cup\partial C^{(2)}_{\blue}$. Since~$|X| \ge 25\eps^{1/12}N > |X\sm X'| + 12\eps^{1/12}N + 12\eps N +\binom{12}{2}\eps N + 13$, by~\eqref{eq:complete-codegree} and Lemma~\ref{lem:setup-properties}\ref{item:dominant-monoxcomp-shadows-nbrhood-sizes} there is a subset~$X''\subseteq X'$ of~$13$ vertices such that~$\binom{X''}{3}\subseteq E(H)$ and~$\binom{X''}{2}\subseteq(\partial C^{(1)}_{\blue}\cup\partial C^{(2)}_{\blue})\cap\partial S$. Now since~$r(C^{(3)}_4)=13$ and~$\binom{X''}{2}\subseteq(\partial C^{(1)}_{\blue}\cup\partial C^{(2)}_{\blue})\cap\partial S$, by~\ref{item:scarlet-no-C4} there is a monochromatic blue copy of~$C^{(3)}_4$ on~$H[X'']$ in~$C^{(1)}_{\blue}$ or~$C^{(2)}_{\blue}$.
\end{proof}

\begin{proof}[Proof of Lemma~\ref{lem:no-dominant-pair-scarlet-C5}]
    We first show that the blue components can be grouped up into three large sets.

    \begin{claim} \label{claim:three-real-parts}
        There is a partition~$\binom{V}{2}=F_1\cup F_2\cup F_3$ such that
        \begin{enumerate}[label=(PP\arabic{*})]
            \item \label{item:partition-class} $F_1$, $F_2$ and~$F_3$ are pairwise disjoint,
            \item \label{item:partition-real} for every~$i\in[3]$ we have $\left|F_i\cap\binom{V'_{\red}}{2}\right|\ge\eps^{1/3}N^2$, and
            \item \label{item:partition-components} for every blue component~$C_{\blue}$ there exists~$i\in[3]$ such that~$\partial C_{\blue}\subseteq F_i$.
        \end{enumerate}
    \end{claim}

    \begin{claimproof}
        We first obtain a partition of~$\binom{V'_{\red}}{2}$ as follows. For each blue component~$C_{\blue}$ let~$\partial C_{\blue}\cap\binom{V'_{\red}}{2}$ be a partition class, and let each pair which is in no blue edge form a partition class by itself; this is well-defined by~\eqref{eq:components-disjoint-shadows}. Enumerate the sizes of the partition classes in decreasing order as~$a_1\ge\dots\ge a_k$. Since~\eqref{eq:blue-monoxcomps-dominant-pair} does not hold for any pair~$C^{(1)}_{\blue}$ and~$C^{(2)}_{\blue}$ of blue components, we have~$a_1+a_2\le\binom{|V'_{\red}|}{2}-2\eps^{1/3}N^2$ and~$k\ge3$. Now~\ref{item:few-blue-vertices} and Lemma~\ref{lem:setup-properties}\ref{item:colour-bipartition-near-complete} yield~$\left|\binom{V'_{\red}}{2}\right| \ge 5\eps^{1/3}N^2$, so by applying Lemma~\ref{lem:parts-allocations}\ref{item:three-parts} we may group up the terms~$a_i$ into three sums which are all at least~$\eps^{1/3}N^2$. Let~$F'_1$,$F'_2$ and~$F'_3$ be the subsets of~$\binom{V'_{\red}}{2}$ corresponding to these sums of terms~$a_i$; clearly, they partition~$\binom{V'_{\red}}{2}$. Finally, arbitrarily extend the partition~$\binom{V'_{\red}}{2} = F'_1\cup F'_2\cup F'_3$ to a partition~$\binom{V}{2}=F_1\cup F_2\cup F_3$ such that for each~$i\in[3]$ we have~$F'_i\subseteq F_i$ and for each blue component~$C_{\blue}$ such that~$\partial C_{\blue}\cap\binom{V'_{\red}}{2}\subseteq F'_i$ we have~$\partial C_{\blue}\subseteq F_i$; since~$F'_i=F_i\cap\binom{V'_{\red}}{2}$ for all~$i\in[3]$, this clearly satisfies the requirements.
    \end{claimproof}

    For each~$i\in[3]$ and~$v\in V'_{\red}$ define~$F'_i:=F_i\cap\binom{V'_{\red}}{2}$ and~$\deg_i(v):=\deg_{F'_i}(v)$.

    \begin{claim} \label{claim:vertex-two-parts-only}
        For each vertex~$v\in V'_{\red}$ there is an index~$i\in[3]$ such that~$\deg_i(v)\le7\eps N$.
    \end{claim}

    \begin{claimproof}
        Suppose that there is a vertex~$v\in V'_{\red}$ such that for each~$i\in[3]$ we have~$\deg_i(v)>7\eps N$. For each~$i\in[3]$ let~$U_i$ denote the set of vertices~$w$ in the neighbourhood of~$v$ in~$F'_i$ such that~$vw\in\partial S$; let~$U:=U_1\cup U_2\cup U_3$. By Lemma~\ref{lem:setup-properties}\ref{item:dominant-monoxcomp-shadows-nbrhood-sizes} we have~$|V'_{\red}\sm U|\le\eps N$, so for all~$i\in[3]$ we have~$|U_i|>6\eps N$.

        We say that a triple~$uu'u''$ is \emph{$v$-good} if~$H[\{v,u,u',u''\}]$ is a copy of~$C^{(3)}_4$. Note that all~$v$-good edges~$u_1u_2u_3\in E(H)$ such that~$u_i\in U_i$ for all~$i\in[3]$ are blue and lie in the same blue component~$C^*_{\blue}$. Indeed, let~$u_1u_2u_3\in E(H)$ be a~$v$-good edge such that for all~$i\in[3]$ we have~$u_i\in U_i$. Note that all the edges~$vu_iu_j$ with~$1\le i<j\le 3$ must be red edges in~$S$ by~\ref{item:partition-components} and the definition of~$U_i$. Since by~\ref{item:scarlet-no-C4} we have~$C^{(3)}_4\nsubseteq S$, the edge~$u_1u_2u_3$ must be blue. Now let~$u'_1u'_2u'_3$ and~$u''_1u''_2u''_3$ be two~$v$-good edges such that for all~$i\in[3]$ we have~$u'_i,u''_i\in U_i$. By~\eqref{eq:complete-codegree} applied to the pairs~$vu'_1,vu'_2,u'_1u'_2,vu''_1,vu''_2,u''_1u''_2$ and because~$|U_3|>6\eps N$, there is a vertex~$x_3\in U_3$ such that both~$u'_1u'_2x_3$ and~$u''_1u''_2x_3$ are~$v$-good edges. Similarly, by~\eqref{eq:complete-codegree} applied to the pairs~$vu'_1,vx_3,u'_1x_3,vu''_1,u''_1x_3$ and because~$|U_2|>6\eps N$, there is a vertex~$x_2\in U_2$ such that both~$u'_1x_2x_3$ and~$u''_1x_2x_3$ are~$v$-good edges. Hence, $u'_1u'_2u'_3$ and~$u''_1u''_2u''_3$ are in the same blue component. 

        Let~$1\le i<j\le 3$, and let~$u_i\in U_i$ and~$u_j\in U_j$ be a pair of vertices such that~$vu_iu_j\in E(H)$. Note that the edge~$vu_iu_j$ must be red and in~$S$; in particular, we have~$u_iu_j\in\partial S$. Taking the unique~$k\in[3]\sm\{i,j\}$, by~\eqref{eq:complete-codegree} and because~$|U_k|>6\eps N$ there exists~$u_k\in U_k$ such that~$u_iu_ju_k\in E(H)$ is~$v$-good. Since~$u_iu_ju_k\in C^*_{\blue}$ by the previous paragraph, we have~$u_iu_j\in\partial C^*_{\blue}$. Hence, the shadow~$\partial C^*_{\blue}$ contains all pairs of vertices~v with~$1\le i<j\le 3$ and~$vu_iu_j\in E(H)$; without loss of generality, we may assume that~$u_iu_j\in F'_1$.

        Let~$i,j\in[3]$ satisfy~$j\ne i$, and let~$u_iu'_i\in\binom{U_i}{2}$ be a pair such that~$vu_iu'_i\in E(H)$. By~\eqref{eq:complete-codegree} applied to the pairs~$vu_i,vu'_i,u_iu'_i$ and because~$|U_j|>3\eps N$, there is a vertex~$u_j\in U_j$ such that~$u_iu'_iu_j\in E(H)$ is~$v$-good. The edges~$vu_iu_j$ and~$vu'_iu_j$ are both in~$S$ and by~\ref{item:scarlet-no-C4} we have~$C^{(3)}_4\nsubseteq S$, so~$vu_iu'_i$ or~$u_iu'_iu_j$ is blue. Since~$vu_i\in F'_i$ and~$u_iu_j\in F'_1$, we have~$u_iu'_i\in F'_i\cup F'_1$.

        Let~$i\in\{2,3\}$. By~\ref{item:partition-real} we have~$|F'_i|\ge\eps^{1/3}N^2$. Since~$|V'_{\red}\sm U|\le\eps N$, there are at most~$\eps N^2$ pairs in~$F'_i\sm\binom{U}{2}$. Now for all~$u\in U$ we have~$vu\in\partial H$, so by~\eqref{eq:complete-codegree} there are at most~$\eps N^2$ pairs~$u_iu'_i\in\binom{U}{2}$ such that~$vu_iu'_i\notin E(H)$. Since the previous two paragraphs implies that every pair~$uu'$ in~$F'_i\cap\binom{U}{2}$ with~$vuu'\in E(H)$ is in~$\binom{U_i}{2}$, we conclude that there are at least~$\frac{1}{2}\eps^{1/3}N^2$ pairs~$u_iu'_i$ of~$F'_i$ contained in~$U_i$ such that~$vu_iu'_i\in E(H)$; in particular, there is a matching~$M_i$ in~$F'_i$ of size~$5\eps N$ such that every~$u_iu'_i\in M_i$ satisfies~$vu_iu'_i\in E(H)$. We claim that there are pairs~$u_2u'_2\in M_2$ and~$u_3u'_3\in M_3$ such that~$H[\{u_2,u'_2,u_3,u'_3\}]$ is a copy of~$C^{(3)}_4$. Indeed, this follows by applying Lemma~\ref{lem:bipartite-high-degree-shared-edge} and noting that for~$i\in\{2,3\}$ by~\eqref{eq:complete-codegree} each~$u_iu'_i\in M_i$ has at most~$\eps N$ pairs~$yy'\in M_{5-i}$ such that~$u_iu'_iy\notin E(H)$ or~$u_iu'_iy'\notin E(H)$.

        By the previous paragraph, there are~$u_2,u'_2\in U_2$ and~$u_3,u'_3\in U_3$ such that~$H[\{u_2,u'_2,u_3,u'_3\}]$ is a copy of~$C^{(3)}_4$ and for~$i\in\{2,3\}$ we have~$vu_iu'_i\in E(H)$ and~$u_iu'_i\in F'_i$. Now by~\ref{item:scarlet-no-C4} we have~$C^{(3)}_4\nsubseteq S$ and we previously established that~$\{u_2u_3,u'_2u_3,u_2u'_3,u'_2u'_3\}\subseteq\partial S$, so at least one of the four edges induced by~$\{u_2,u'_2,u_3,u'_3\}$, say~$u_2u'_2u_3$, must be blue. Since~$u_2u_3\in F'_1$, we have~$u_2u'_2\in F'_1\cap F'_2$, which contradicts~\ref{item:partition-class}.
    \end{claimproof}

    For~$1\le i<j\le 3$ set~$W_{ij}:=\{v\in V'_{\red}:\deg_i(v),\deg_j(v)>7\eps N\}$; note that Claim~\ref{claim:vertex-two-parts-only} implies that~$W_{12}$, $W_{13}$ and~$W_{23}$ are pairwise disjoint.

    \begin{claim} \label{claim:pair-part-robust}
        For every~$1\le i<j\le 3$ we have~$|W_{ij}|\ge24\eps^{1/2}N$.
    \end{claim}

    \begin{claimproof}
        Since the arguments for the other cases are analogous by symmetry, we restrict ourselves to the case~$i=1$ and~$j=2$. For~$\ell\in[2]$ define the graph~$G_\ell=(V'_{\red},F'_\ell\cap\partial S)$; it has at least~$\frac{1}{2}\eps^{1/2}N^2$ edges by~\ref{item:partition-real} and Lemma~\ref{lem:setup-properties}\ref{item:dominant-monoxcomp-shadows-nbrhood-sizes}. By iteratively removing vertices of low degree, for each~$\ell\in[2]$ we obtain a subgraph~$G'_\ell\subseteq G_\ell$ with minimum degree at least~$\eps^{1/3}N^2/(4N)\ge1152\eps^{1/2}N$ and greedily pick a matching~$M'_\ell$ of size at least~$576\eps^{1/2}N$ in~$G'_\ell$. Observe that we may select submatchings~$M_1\subseteq M'_1$ and~$M_2\subseteq M'_2$ of size~$m:=144\eps^{1/2}N$ such that~$M_1\cup M_2$ is a matching. Indeed, pick~$M_1\subseteq M'_1$ and~$M_2\subseteq M'_2$ with maximum~$|M_1|=|M_2|=m'$ and note that if~$m'<m$ then we could enlarge both~$M_1$ and~$M_2$ because~$M_1\cup M_2$ would intersect at most~$4m'<4m=|M'_1|=|M'_2|$ edges of each of~$M'_1$ and~$M'_2$.

        For each~$vv'\in M_\ell$ with~$\ell\in[2]$, by~\eqref{eq:complete-codegree} there are at most~$2\eps N$ pairs~$ww'\in M_{3-\ell}$ such that~$vv'w\notin E(H)$ or~$vv'w'\notin E(H)$. Hence, there are at least~$m(m-2\eps N)\ge 2m^2/3$ pairs~$v_1v'_1\in M_1$ and~$v_2v'_2\in M_2$ such that~$H[\{v_1,v'_1,v_2,v'_2\}]$ is a copy of~$C^{(3)}_4$. Since by~\ref{item:scarlet-no-C4} we have~$C^{(3)}_4\nsubseteq S$, one of the edges must be blue and so two of the pairs~$v_1v_2,v_1v'_2,v'_1v_2,v'_1v'_2$ belong to the same shadow of a blue component from~$F'_1\cup F'_2$. Hence, there exists~$\ell\in[2]$ such that the induced bipartite subgraph of~$G'_\ell$ on~$V(M_1)\cup V(M_2)$ has at least~$2m^2/3$ edges. By iteratively removing vertices of low degree, there are at least~$m/6\ge24\eps^{1/2}N$ vertices of~$V(M_{3-\ell})$ with minimum degree at least~$m/6>7\eps N$ in~$G'_\ell$. Since the vertices of~$V(M_{3-\ell})$ have degree at least~$1152\eps^{1/2}N$ in~$G'_{3-\ell}$, we have~$|W_{12}|\ge24\eps^{1/2}N$.
    \end{claimproof}

    We first claim that at most~$14\eps N^2$ pairs between~$W_{ij}$ and~$W_{i'j'}$ are not from~$F'_{\ell}$, where~$\{\ell\}=\{i,j\}\cap\{i',j'\}$; since the other cases are analogous by symmetry, we restrict ourselves to the case~$\ell=1$. Indeed, suppose that more than~$14\eps N^2$ pairs between~$W_{12}$ and~$W_{13}$  are not from~$F'_1$. Then more than~$7\eps N^2$ pairs come from, say, $F'_3$. Hence, there exists a vertex~$v\in W_{12}$ such that~$\deg_3(v)>7\eps N$. But this contradicts Claim~\ref{claim:vertex-two-parts-only}.

    Next, we shall pick distinct vertices~$w_{12},w'_{12}\in W_{12}$, $w_{13},w'_{13}\in W_{13}$, $w_{23},w'_{23}\in W_{23}$ such that, setting~$Y:=\{w_{12},w'_{12},w_{13},w'_{13},w_{23},w'_{23}\}$, the following hold.
    \begin{enumerate}[label=(F\arabic{*})]
        \item \label{item:six-clique} $H[Y]$ is a copy of~$K^{(3)}_6$.
        \item \label{item:azure-pair-clique} $\binom{Y}{2}\subseteq \partial S$.
        \item \label{item:red-cross-pairs} For each~$\ell\in[3]$ and distinct pairs~$\{i,j\},\{i',j'\}\in\binom{[3]}{2}$ such that~$\{i,j\}\cap\{i',j'\}=\{\ell\}$ we have~$w_{ij}w_{i'j'},w_{ij}w'_{i'j'},w'_{ij}w_{i'j'},w'_{ij}w'_{i'j'}\in F'_\ell$.
    \end{enumerate}
    Indeed, the previous paragraph combined with Claim~\ref{claim:pair-part-robust}, Lemma~\ref{lem:setup-properties}\ref{item:dominant-monoxcomp-shadows-nbrhood-sizes} and~\eqref{eq:complete-codegree} implies that there are distinct vertices~$w_{12},w'_{12}\in W_{12}$ and subsets~$W'_{13}\subseteq W_{13}$ and~$W'_{23}\subseteq W_{23}$ such that~$w_{12}w'_{12}\in\partial S$, $|W'_{13}|,|W'_{23}|\ge22\eps^{1/2}N$ and for all~$i\in\{1,2\}$ and~$x\in W'_{i3}$ we have~$w_{12}w'_{12}x\in E(H)$ and~$w_{12}x,w'_{12}x\in\partial S\cap F'_i$. Now by the previous paragraph, Lemma~\ref{lem:setup-properties}\ref{item:dominant-monoxcomp-shadows-nbrhood-sizes} and~\eqref{eq:complete-codegree}, there are distinct vertices~$w_{13},w'_{13}\in W'_{13}$ and a subset~$W''_{23}\subseteq W'_{23}$ such that~$w_{13}w'_{13}\in\partial S$, $w_{12}w_{13}w'_{13},w'_{12}w_{13}w'_{13}\in E(H)$, $|W''_{23}|\ge20\eps^{1/2}N$ and for all~$x\in W''_{23}$ we have~$w_{13}x,w'_{13}x\in\partial S\cap F'_3$ and~$w_{12}w_{13}x,w_{12}w'_{13}x,w'_{12}w_{13}x,w'_{12}w'_{13}x,w_{13}w'_{13}x\in E(H)$. Finally, by Lemma~\ref{lem:setup-properties}\ref{item:dominant-monoxcomp-shadows-nbrhood-sizes} and~\eqref{eq:complete-codegree}, there are vertices~$w_{23},w'_{23}\in W''_{23}$ such that~$w_{23}w'_{23}\in\partial S$ and for all~$i\in\{2,3\}$ we have $w_{1i}w_{23}w'_{23},w'_{1i}w_{23}w'_{23}\in E(H)$. Hence, we have a set~$Y$ of vertices which satisfies~\ref{item:six-clique}--\ref{item:red-cross-pairs}.

    To finish the proof, we show that~\ref{item:scarlet-no-C4} and~\ref{item:six-clique}--\ref{item:red-cross-pairs} implies that~$H[Y]$ contains a copy of~$C^{(3)}_5$ in~$S$. Without loss of generality, we may assume that~$w_{12}w'_{12}\notin F'_1$ because~$F'_1\cap F'_2=\nth$. Since by~\ref{item:red-cross-pairs}~$w_{12}w_{13},w_{12}w'_{13}\in F'_1$, the edges~$w_{12}w'_{12}w_{13}$ and~$w_{12}w'_{12}w'_{13}$ are red. Now by~\ref{item:scarlet-no-C4}, \ref{item:six-clique} and~\ref{item:azure-pair-clique} there is a blue edge in~$H[w_{12},w'_{12},w_{13},w'_{13}]$, so at least one of~$w_{12}w_{13}w'_{13}$ and~$w'_{12}w_{13}w'_{13}$ is blue. Since by~\ref{item:red-cross-pairs}~$w_{12}w_{13},w'_{12}w_{13}\in F'_1$, we have~$w_{13}w'_{13}\in F'_1$. Analogously, we have~$w_{23}w'_{23}\in F'_3$ and~$w_{12}w'_{12}\in F'_2$.

    Observe that by~\ref{item:red-cross-pairs} all edges of the form~$x_{12}x_{13}x_{23}$, where for all~$1\le i<j\le 3$ we have~$x_{ij}\in\{w_{ij},w'_{ij}\}$, are red because all pairs in a blue edge must belong to the same~$F'_i$; by the same reasoning, the edges~$w_{12}w'_{12}w_{13}$, $w_{12}w'_{12}w'_{13}$, $w_{13}w'_{13}w_{23}$, $w_{13}w'_{13}w'_{23}$, $w_{23}w'_{23}w_{12}$ and~$w_{23}w'_{23}w'_{12}$ are all red. In particular, $H[Y]$ contains a red copy of~$C^{(3)}_5$ in~$S$ with the edges~$w_{12}w_{23}w_{13}$, $w_{23}w_{13}w'_{13}$, $w_{13}w'_{13}w'_{23}$, $w'_{13}w'_{23}w_{12}$ and~$w'_{23}w_{12}w_{23}$.
\end{proof}

\end{document}